\documentclass[reqno]{amsart}
\usepackage{amsmath, colortbl}
\usepackage{amsfonts}
\usepackage{amssymb}
\usepackage{mathrsfs}
\usepackage{graphicx}
\usepackage{flafter}
\newtheorem{thm}{Theorem}[section]
\newtheorem{cor}[thm]{Corollary}
\newtheorem{lem}[thm]{Lemma}

\theoremstyle{definition}

\newtheorem{rem}[thm]{Remark}
\numberwithin{equation}{section}

\begin{document}
\Large{
\thispagestyle{empty}

\vspace{1 true cm}

\title[   Regularity of the solution of the  3D Navier-Stokes
Equations ] {Some new regularity criteria \\ for the  3D Navier-Stokes
Equations}%
\author[Daoyuan Fang \ Chenyin Qian]{Daoyuan Fang \hspace{0.5cm}Chenyin Qian\\
Department of Mathematics, Zhejiang University,\\ Hangzhou, 310027,
China}

\thanks{$^*$E-mail addresses: dyf@zju.edu.cn (D. Fang), qcyjcsx@163.com(C. Qian)}

\begin{abstract}
Several types of  new regularity criteria  of Leray-Hopf weak
solutions $u$ to the 3D Navier-Stokes equations are obtained.  Some
of them are based on the third component $u_3$ of velocity under the
Prodi-Serrin index condition.  And a very recent work of the
authors, based on only one of the nine entries of the gradient
tensor, is renovated. At last,  some regularity criteria which are
dependent on some parameter $\epsilon$ are obtained.
\end{abstract}

\maketitle

\section{Introduction}
 In  the present  paper, we  address sufficient  conditions for the regularity of  weak solutions of the Cauchy problem for the
Navier-Stokes equations in $ \mathbb{R}^{3}\times (0,T)$:
\begin{equation} \label{a}
 \left\{\begin{array}{l}
\displaystyle \frac{\partial u}{\partial t}-\nu \Delta
u+(u\cdot\nabla)u+\nabla p=0,
\ \mbox{\ in}\ \mathbb{R}^{3}\times (0,T),\vspace{1mm}\\
\displaystyle\nabla\cdot u=0,\hspace{3.52cm} \mbox{\ in}\ \mathbb{R}^{3}\times (0,T),\\
\displaystyle u(x, 0)=u_{0},\hspace{0.2cm}\mbox{\ in}\ \mathbb{R}^{3},\\
\end{array}
\right. \end{equation} where $u=(u_{1},u_{2},u_{3}):
\mathbb{R}^{3}\times (0,T)\rightarrow \mathbb{R}^{3}$ is the
velocity field, $ p: \mathbb{R}^{3}\times (0,T)\rightarrow
\mathbb{R}^{3}$ is a scalar pressure, and $u_{0}$ is the initial
velocity field, $\nu>0$ is the viscosity.  We set
$\nabla_{h}=(\partial_{x_{1}},\partial_{x_{2}})$ as the horizontal
gradient operator and
$\Delta_{h}=\partial_{x_{1}}^{2}+\partial_{x_{2}}^{2}$ as the
horizontal Laplacian, and $\Delta$ and $\nabla$ are the usual
Laplacian and the gradient operators, respectively.  Here we use the
classical notations
$$
(u\cdot\nabla)v=\sum_{i=1}^{3}u_{i}\partial_{x_{i}}v_{k}, \ (
k=1,2,3),\ \ \ \nabla\cdot u=\sum_{i=1}^{3}\partial_{x_{i}}u_{i},
$$
 and for sake of simplicity,  we denote $\partial_{x_{i}}$
by $\partial_{{i}}$.

\par   It is well known that the weak solution of the Navier-Stokes equations \eqref{a} is  unique  and
regular in two dimensions. However, in three dimensions, the
regularity problem of weak solutions of Navier-Stokes equations is
an outstanding open problem in mathematical fluid mechanics. The
weak solutions are known to exist globally in time, but the
uniqueness, regularity, and continuous dependence on initial data
for weak solutions are still open problems. Furthermore, strong
solutions in the 3D case are known to exist for a short interval of
time whose length depends on the initial data. Moreover, this strong
solution is known to be unique and to depend continuously on the
initial data (see, for example, \cite{[21]}, \cite{[22]}). Let us
recall the definition of Leray-Hopf weak solution. We set
$$
\mathcal {V}=\{\phi: \ \mbox{ the 3D vector valued}\ C_{0}^{\infty}
\ \mbox{functions and}\ \nabla\cdot\phi=0\},
$$
which will form the space of test functions. Let $H$ and $V$ be the
closure spaces of $\mathcal {V}$ in $L^2$ under $L^2$-topology, and
in $H^1$ under $H^1$-topology, respectively.
\par
For $u_0\in H$,  the existence of weak solutions of \eqref{a} was
established by Leray  \cite{[14]} and Hopf in  \cite{[9]}, that is,
$u$
satisfies the following properties:\\
(i) $u\in C_{w}([0,T); H)\cap L^{2}(0,T; V)$, and $\partial_{t}u\in
L^{1}(0,T; V^{\prime})$, where $V^{\prime}$ is the dual space of
$V$;\\
(ii) $u$ verifies (1.1) in the sense of distribution, i.e., for
every test function $\phi\in C_{0}^{\infty}([0,T);\mathcal {V})$,
and for almost every $t, t_{0}\in (0,T)$, we have
$$
\begin{array}{ll}
 \displaystyle&\displaystyle\int_{\mathbb{R}^{3}} u(x,t)\cdot\phi(x,t)dx-\int_{\mathbb{R}^{3}}
u(x,t_{0})\cdot\phi(x,t_{0})dx\vspace{2mm}\\
\displaystyle &\ \ \ \
=\displaystyle\int_{t_{0}}^{t}\int_{\mathbb{R}_{3}}[u(x,t
)\cdot(\phi_{t}(x,t)+\nu\Delta\phi(x,t))]dxds\vspace{2mm}\\
&\ \ \ \ \ \ \
+\displaystyle\int_{t_{0}}^{t}\int_{\mathbb{R}_{3}}[(u(x,t)\cdot\nabla)\phi(x,t)]\cdot
u(x,t))]dxds
\end{array}
$$
 (iii) The energy inequality,
i.e.,
$$
\|u(\cdot,t)\|_{L^{2}}^{2}+2\nu\int_{t_0}^{t}\|\nabla
u(\cdot,s)\|_{L^{2}}^{2}ds\leq\|u_{0}\|_{L^{2}}^{2},
$$
for every $t$ and almost every $t_{0}$.

It is well known, if $u_{0}\in V$, a weak solution becomes  strong
solution of (1.1) on $(0, T)$ if, in addition, it satisfies
$$
u\in C([0,T); V)\cap L^{2}(0,T; H^{2}) \ \mbox{and}\
\partial_{t}u\in L^{2}(0,T; H).
$$
We know the strong solution is regular(say, classical) and unique
(see, for example, \cite{[21]}, \cite{[22]}).
\par  Researchers are interested in the classical problem of finding sufficient conditions for
weak solutions of (1.1) such that the weak solutions become regular,
and the first result is usually referred as Prodi-Serrin conditions
(see \cite{[18]} and \cite{[20]}), which states that if a weak
solution $u$ is in the class of
\begin{equation}\label{19d} u\in L^{t}(0,T; L^{s}(\mathbb{R}^{3})),\ \
\frac{2}{t}+\frac{3}{s}=1,\ s\in[3,\infty],
\end{equation} then the weak solution becomes regular.
Recently,  H. Bae and  H. Choe in \cite{[207]} gave a two components
Prodi-Serrin index criterion. Up to now, there are many results show that
one can use only one component (say $u_3$) to determine the
regularity of $u$. Say, I. Kukavica and M. Ziane in \cite{[10]}
proved a regularity criterion  under the following condition
$$ u_3\in L^{t}(0,T; L^{s}(\mathbb{R}^{3})),\ \
\frac{2}{t}+\frac{3}{s}\leq \frac{5}{8},\ s\in[\frac{24}{5},\infty].
$$
Then, it was improved by C. Cao and  E.Titi in
\cite{[3]}  to
$$ u_3\in L^{t}(0,T; L^{s}(\mathbb{R}^{3})),\ \
\frac{2}{t}+\frac{3}{s}\leq \frac{2(s+1)}{3s},\
s\in(\frac{7}{2},\infty].
$$
And then,  Y. Zhou and M. Pokorn$\acute{\mbox{y} }$ in \cite{[26]}
changed the regular criterion to
$$ u_3\in L^{t}(0,T; L^{s}(\mathbb{R}^{3})),\ \
\frac{2}{t}+\frac{3}{s}\leq \frac{3}{4}+\frac{1}{2s},\
s\in(\frac{10}{3},\infty].
$$
More relative results, we refer to \cite{[16]},
\cite{[25]} and the reference there in. One can see  that the above
mentioned results on $u_3$ cannot satisfy the Prodi-Serrin index
condition, and it seems to be a price when one reduce the components
of $u$ to one. It is nature to think about what supplement is
necessary to insure the Prodi-Serrin  condition based on one
velocity component. For example, P. Penel and M.
Pokorn$\acute{\mbox{y} }$ in \cite{[37]} proved the $u$ was regular,
if
$$
u_3\in L^{\beta}(0,T; L^{\alpha}(\mathbb{R}^3))\ \mbox{with}\
\frac{3}{\alpha}+\frac{2}{\beta}\leq 1, \alpha\in(3,\infty],
\beta\in [2,\infty),
$$
 and one of the following conditions holds true:\\
$\hspace{2cm}(a) \ \partial_3u_1, \partial_3u_2$ belong to
$L^{\beta}(0,T; L^{\alpha}(\mathbb{R}^3))\ \mbox{with}\
\frac{3}{\alpha}+\frac{2}{\beta}\leq 2, \alpha\in(3/2,\infty],
\beta\in [1,\infty)$;\\
$\hspace{2cm}(b) \ \partial_2u_1, \partial_1u_2$ belong to
$L^{\beta}(0,T; L^{\alpha}(\mathbb{R}^3))\ \mbox{with}\
\frac{3}{\alpha}+\frac{2}{\beta}\leq 2, \alpha\in[2,3],
\beta\in [2,\infty]$;\\
$\hspace{2cm}(c) \ \partial_3u_2$ belong to $L^{\beta}(0,T;
L^{\alpha}(\mathbb{R}^3))\ \mbox{with}\
\frac{3}{\alpha}+\frac{2}{\beta}\leq 2, \alpha\in(3/2,\infty],
\beta\in [1,\infty)$, and  $\ \partial_2u_1$ belong to
$L^{\beta}(0,T; L^{\alpha}(\mathbb{R}^3))\ \mbox{with}\
\frac{3}{\alpha}+\frac{2}{\beta}\leq 2, \alpha\in[2,3],
\beta\in [2,\infty].$\\
Moreover, the authors also mentioned in the Remark 2  in
\cite{[37]}, the condition in $(a)$ can be replaced by $
\partial_3u_2$, $\ \partial_2u_2$, or  $
\partial_3u_2$, $\ \partial_1u_1$, or  $
\partial_3u_1$, $\ \partial_2u_2$, or  $
\partial_3u_1$, $\ \partial_2u_1$. Similarly, in $(c)$ one
can replace  $
\partial_3u_2$ by $\ \partial_3u_1$, and replace
 $
\partial_2u_1$ by $\ \partial_1u_2$ respectively.

\par From the above
, we can see  that the assumptions on derivative component did not
contain $\partial_iu_3\ i=1,2,3$. One purpose of this paper is to
capture this work, by using the incompressibility condition. We give
an estimate on velocity, which is different from \cite{[37]} and
then get a regularity criterion  on $\partial_iu_3\ i=1,2,3$, for
detail see the proof of Theorem \ref{t1.3} below. On the other hand,
similar to $(a)$, $(b)$ and $(c)$,  we also consider cases of  the
given conditions in terms of  only one component $\partial_iu_j$ of
$\nabla u$ such that $u_3$ satisfies the Prodi-Serrin condition in
Theorem \ref{t1.4} and Corollary \ref{c1.1}.
\par
Besides, we  would like to point out that the full regularity of
weak solutions can also be proved under alternative assumptions on
the gradient of the velocity $\nabla u$, for instance
\begin{equation}\label{9d}
\nabla u\in L^{t}(0,T; L^{s}(\mathbb{R}^{3})),\ \
\frac{2}{t}+\frac{3}{s}=2,\ s\in[\frac{3}{2},\infty].
\end{equation}
Enlightened by the above, we also want to get some better regularity
criteria which are also  coincident with the standard Prodi-Serrin
condition based on some components of $\nabla u$. To begin with, we mention some results in this direction at first, P. Penel and M.
Pokorn$\acute{\mbox{y} }$ in \cite{[37]} proved that if
$$
\partial_3u\in L^{\beta}(0,T;
L^{\alpha}(\mathbb{R}^3))\ \mbox{with}\
\frac{3}{\alpha}+\frac{2}{\beta}\leq \frac{3}{2},
\alpha\in[2,\infty], \beta\in [1,\infty).
$$
then the weak solution was regular. After that many authors improved
this result, such as I. Kukavica and M. Ziane in  \cite{[1023]}
considered the case of the condition
$$
\partial_3u\in L^{\beta}(0,T;
L^{\alpha}(\mathbb{R}^3))\ \mbox{with}\
\frac{3}{\alpha}+\frac{2}{\beta}\leq 2, \alpha\in[\frac{9}{4},3].
$$
As to the gradient of one velocity component $\nabla u_3$, M.
Pokorn$\acute{\mbox{y} }$ in \cite{[17]} proved the weak solution
was actually regular if  $\nabla u_3$ satisfied
$$
\nabla u_3\in L^{\beta}(0,T; L^{\alpha}(\mathbb{R}^3))\ \mbox{with}\
\frac{3}{\alpha}+\frac{2}{\beta}\leq \frac{3}{2},
\alpha\in[2,\infty].
$$
 Y. Zhou and M. Pokorn$\acute{\mbox{y} }$ in \cite{[26]}
improved the result to
$$
\nabla u_3\in L^{\beta}(0,T; L^{\alpha}(\mathbb{R}^3))\ \mbox{with}\
\frac{3}{\alpha}+\frac{2}{\beta}\leq  \left\{\begin{array}{l}
\displaystyle \frac{19}{12}+\frac{1}{2\alpha}, \alpha\in(\frac{30}{19},3]\vspace{2mm}\\
\displaystyle\frac{3}{2}+\frac{3}{4\alpha}, \alpha\in(3,\infty],\\
\end{array}
\right.
$$
moreover, Y. Zhou and M. Pokorn$\acute{\mbox{y} }$ also proved a
improved result, for more  detail we refer to \cite{[27]}. Motivated
by the above,  we  consider the case of  two gradient velocity components and one of them satisfies
the Prodi-Serrin condition, see Theorem \ref{t1.9} and
Corollary \ref{c1.2}. We shall point out that two gradient velocity
components are not all the diagonal elements, this is more
difficult than the  diagonal case, for the detail see Remark
\ref{r1.3} below.
\par
 In  \cite{[26]}, the authors also studied the
regularity of the solutions of the Navier-Stokes equations under the
assumption on $\partial_{3}u_{3}$, namely,
\begin{equation}\label{jbb}
\partial_{3}u_{3}\in L^{\beta}(0,T; L^{\alpha}(\mathbb{R}^{3})),\
\frac{3}{\alpha}+\frac{2}{\beta}< \frac{4}{5},\
\alpha\in(\frac{15}{4}, \infty].
\end{equation} Recently, the regularity criterion in terms of only one of the gradient tensor was gotten by C. Cao and E.
 Titi in \cite{[2]} under the  assumptions
$$
\frac{\partial u_{j}}{\partial x_{k}}\in L^{\beta}(0,T;
L^{\alpha}(\mathbb{R}^{3})), \ \mbox{when} \ j\neq k
$$
\begin{equation}\label{c}
\mbox{and where}  \ \alpha>3, 1\leq\beta<\infty,\  \mbox{and}\ \
\frac{3}{\alpha}+\frac{2}{\beta}\leq\frac{\alpha+3}{2\alpha},\end{equation}
or
$$
\frac{\partial u_{j}}{\partial x_{j}}\in L^{\beta}(0,T;
L^{\alpha}(\mathbb{R}^{3})),
$$
\begin{equation}\label{d}
\mbox{and where}  \ \alpha>2, 1\leq\beta<\infty,\  \mbox{and} \ \
\frac{3}{\alpha}+\frac{2}{\beta}\leq\frac{3(\alpha+2)}{4\alpha}.\end{equation}
In \cite{[201]}, we improved this result. And here,  we again study
it and  get an improvement of the results of \cite{[201]}, which is
shown in Theorem \ref{t1.5}.  However, it is also noted that the
above conditions are not coincident with the
 Prodi-Serrin condition.

Now,  we list our main results as follows:
\begin{thm}\label{t1.3}
 Let $u$ be a Leray-Hopf weak solution to the 3D Navier-Stokes
equations \eqref{a} with the initial value $u_{0}\in V$. Suppose
 \begin{eqnarray}\label{5} u_3\in L^{\beta_1}(0,T;
L^{\alpha_1}(\mathbb{R}^3))\ \mbox{with}\
\frac{3}{\alpha_1}+\frac{2}{\beta_1}\leq1, \alpha_1\in(3,\infty],
\end{eqnarray}
and one of the following conditions holds:
\begin{eqnarray}\label{6} (i)\  \partial_3u_2,\partial_3u_3\in L^{\beta_2}(0,T;
L^{\alpha_2}(\mathbb{R}^3))\ \mbox{with}\
\frac{3}{\alpha_2}+\frac{2}{\beta_2}\leq2, \alpha_2\in[2,3]
\end{eqnarray}
\begin{eqnarray}\label{6} (ii)\  \partial_3u_1,\partial_3u_3\in L^{\beta_2}(0,T;
L^{\alpha_2}(\mathbb{R}^3))\ \mbox{with}\
\frac{3}{\alpha_2}+\frac{2}{\beta_2}\leq2, \alpha_2\in[2,3]
\end{eqnarray}
Then $u$ is regular.
\end{thm}
\begin{thm}\label{t1.4}
Let $u$ be a Leray-Hopf weak solution to the 3D Navier-Stokes
equations \eqref{a}. Suppose that, for some $i,j $ with $1\leq i\leq
3$ and $1\leq j\leq 2$, $u$ satisfies one of the following
conditions:\par (a) $i\neq j$,  suppose the initial value $u_{0}\in
V\bigcap L^{q}(\mathbb{R}^3)$, where $1<q\leq2$, the solution $u$
satisfies
\begin{eqnarray}\label{15} u_3\in L^{\beta_1}(0,T;
L^{\alpha_1}(\mathbb{R}^3))\ \mbox{with}\
\frac{3}{\alpha_1}+\frac{2}{\beta_1}\leq1, \alpha_1\in(3,\infty],
\end{eqnarray}
and
\begin{eqnarray}\label{16} \partial_iu_j\in L^{\beta_2}(0,T;
L^{\alpha_2}(\mathbb{R}^3)),
\end{eqnarray}
with
\begin{eqnarray}\label{17}
\frac{3}{\alpha_2}+\frac{2}{\beta_2}\leq \frac{3-q}{\alpha_2}+q-1, \
\alpha_2\in(\frac{q}{q-1},\infty].
\end{eqnarray}
\par(b) $i=j$,  suppose the initial value $u_{0}\in V$ and $ u_3$
satisfies the condition \eqref{15}, and
\begin{eqnarray}\label{8} \partial_ju_j\in L^{\beta_3}(0,T;
L^{\alpha_3}(\mathbb{R}^3)),
\end{eqnarray}
with
\begin{eqnarray}\label{9}
\frac{3}{2\alpha_3}+\frac{2}{\beta_3}\leq f(\alpha_3), \
\alpha_3\in(\frac{9}{5},\infty], 
\end{eqnarray}
where
$$
f(\alpha_3)=\frac{\sqrt{24\alpha^2_3-24\alpha_3+9}-2\alpha_3}{2\alpha_3}.
$$
Then $u$ is regular.
\end{thm}
\begin{rem}
When we announced the first version of this article on the arXiv.org, we
 were informed by the authors of \cite{[h7]} that they  finished  the same result as the part $(a)$ of Theorem ~\ref{t1.4} with $q=2$. The above is the improved result with  a parameter $q$ satisfying $1<q\leq 2$. For $q>2,$ in fact, one also can get
some results, for example the
condition \eqref{17} can be replaced by
\begin{eqnarray}\label{sw17}
\frac{3}{4\alpha_2}+\frac{2}{\beta_2}\leq
\frac{\sqrt{52\alpha_2^2-60\alpha_2+9}-4\alpha_2}{4\alpha_2}, \
\alpha_2\in(\frac{3+\sqrt 7}{2},\infty].
\end{eqnarray}
\end{rem}

\begin{cor}\label{c1.1} Suppose that $u_{0}\in V$, and
$u$ is a Leray-Hopf weak solution to the 3D Navier-Stokes equations
\eqref{a}.
 Suppose that, for some $i,j
$ with $1\leq i\leq 3$ and $1\leq j\leq 2$, $u$ satisfies one of the
following conditions:\par (a) $i\neq j$,  $u_3$ satisfies the
condition \eqref{15} and
\begin{eqnarray}\label{124} \partial_iu_j\in L^{\beta_2}(0,T;
L^{\alpha_2}(\mathbb{R}^3)),\
\frac{3}{\alpha_2}+\frac{2}{\beta_2}\leq 1, \ \alpha_2\in(2,\infty]
.
\end{eqnarray}
\par(b) $i=j$,  $ u_3$ satisfies the condition \eqref{15}, and
\begin{eqnarray}\label{123} \partial_ju_j\in L^{\beta_3}(0,T;
L^{\alpha_3}(\mathbb{R}^3)),
\frac{3}{\alpha_3}+\frac{2}{\beta_3}\leq \frac{3}{2}, \
\alpha_3\in[2,6]
.
\end{eqnarray}
Then $u$ is regular.
\end{cor}

If we substitute the condition on $u_3$ by the component of the
gradient of the velocity, we have the following regularity
criterion, which is a further improvement of the above mentioned
results of \cite{[37]}.

\begin{thm}\label{t1.9}
Let $u$ be a Leray-Hopf weak solution to the 3D Navier-Stokes
equations \eqref{a} with the initial value $u_{0}\in V$.   Suppose
\begin{eqnarray}\label{o6} \partial_3u_i\in L^{\beta_1}(0,T;
L^{\alpha_1}(\mathbb{R}^3))\ \mbox{with}\
\frac{3}{\alpha_1}+\frac{2}{\beta_1}\leq 2, \alpha_1\in[2,3],
i=1\  \mbox{or} \ 2,
\end{eqnarray}
 and
\begin{eqnarray}\label{o8} \partial_3u_3\in L^{\beta_2}(0,T;
L^{\alpha_2}(\mathbb{R}^3)),
\end{eqnarray}
with
\begin{eqnarray}\label{o9}
\frac{3}{2\alpha_2}+\frac{2}{\beta_2}\leq f(\alpha_2), \
\alpha_2\in[2,3], 
\end{eqnarray}
where
$$
f(\alpha_2)=\frac{\sqrt{24\alpha_2^2-24\alpha_2+9}-2\alpha_2}{2\alpha_2}.
$$
Then $u$ is regular.
\end{thm}

\begin{cor}\label{c1.2} Suppose that $u_{0}\in V$, and
$u$ is a Leray-Hopf weak solution to the 3D Navier-Stokes equations
\eqref{a}. Assume
\begin{eqnarray}\label{oo6} \partial_3u_i\in L^{\beta_1}(0,T;
L^{\alpha_1}(\mathbb{R}^3))\ \mbox{with}\
\frac{3}{\alpha_1}+\frac{2}{\beta_1}\leq 2, \alpha_1\in[2,3],
i=1\  \mbox{or} \ 2,
\end{eqnarray}
and
\begin{eqnarray}\label{ooo6} \partial_3u_3\in L^{\beta_2}(0,T;
L^{\alpha_2}(\mathbb{R}^3))\ \mbox{with}\
\frac{3}{\alpha_2}+\frac{2}{\beta_2}\leq \frac{3}{2},
\alpha_2\in[2,3]
.
\end{eqnarray}
 Then $u$ is regular.
\end{cor}
\begin{rem} \label{r1.3}
Here we only need
two components of the gradient of the velocity and one of them is not on
the diagonal elements of $\nabla u.$ On the case of the diagonal
elements of $\nabla u$, P. Penel and M. Pokorn$\acute{\mbox{y} }$ in
\cite{[37]} proved the $u$ is regular when
$$\partial_2u_2, \partial_3u_3\in L^{\beta}(0,T;
L^{\alpha}(\mathbb{R}^3)),
$$
and $\alpha,\beta$ satisfied
$$\ \frac{3}{\alpha}+\frac{2}{\beta}\leq2,
\alpha\in(\frac{3}{2},\infty], \beta\in [2,\infty).
$$
Moreover, the condition on $\partial_3u_i$ satisfies the
Prodi-Serrin condition, which is an improvement of the result of P.
Penel and M. Pokorn$\acute{\mbox{y} }$ in \cite{[37]}. Finally, we
note that $\partial_3u_i,$ $i=1$ or $2$, is not the diagonal element
of $\nabla u.$ Thus, we cannot use the method of by multiplying
$u_i$ to the $i$th equation of \eqref{a} to get  the form
$\partial_3u_i,$ $i=1$ or $2$. Therefore, it is more difficult to
get the regularity criterion based on $\partial_iu_i$ and
$\partial_ju_j$, $i,j\in\{1,2,3\}$ with $i\neq j$.
\end{rem}
\begin{thm}\label{t1.5} Let  $u_{0}$ and  $u$ be as   in
Theorem \ref{t1.3}.  Suppose in addition
\begin{eqnarray}\label{13} \partial_ku_k\in L^{\beta}(0,T;
L^{\alpha}(\mathbb{R}^3)), \ \mbox{for some}\
k\in\{1,2,3\}\end{eqnarray} with
\begin{equation} \label{14}
 \frac{3}{2\alpha}+\frac{2}{\beta}\leq
\displaystyle g(\alpha), \ \alpha\in(\frac{9}{5},\infty],
\end{equation} where
$$
g(\alpha)=\frac{\sqrt{289\alpha^2-264\alpha+144}-7\alpha}{8\alpha}.
$$
 Then $u$ is
regular.
\end{thm}
\begin{rem} \label{r1.4}

This theorem is an improvement  of \cite{[2]} and \cite{[26]} (see
figure 1 below), and is also an improvement of  Theorem 1.2 (i) and
Theorem 1.3 in \cite{[201]}. Moreover, we point out that the first
part of Theorem 1.2 in
\cite{[201]}, can be simplified to the following form:\\
For $ j\neq k$,  suppose that $u$ satisfies
\begin{eqnarray}\label{z} \int_{0}^{T}\|
\partial_{j}u_{k}\|_{\alpha}^{\beta}d\tau\leq M, \ \mbox{for some} \ M>0, \end{eqnarray}
with
\begin{equation} \label{f}
 \frac{3}{2\alpha}+\frac{2}{\beta}\leq
\displaystyle f(\alpha), \  \alpha\in (3,\infty)\ \mbox {and}\
1\leq\beta<\infty,  \end{equation} where
\begin{equation}\label{m}
f(\alpha)=\frac{\sqrt{103\alpha^2-12\alpha+9}-9\alpha}{2\alpha},
\end{equation}
then  $u$ is  regular. This function shows the same line as in the
Figure 1. in \cite{[201]}. In the proof of Theorem 1.2 (i) of
\cite{[201]}, if we substitute $\sigma_1$ and $\beta$ by
$$\frac{1}{\sigma_1}=\frac{(7-\frac{3}{\alpha})+\sqrt{\frac{9}{\alpha^2}-\frac{12}{\alpha}+103}}{18},
\ \beta=\frac{2\sigma_1}{9-8\sigma_1},$$ we can get the desired
result.
\end{rem}
The following theorems show the variation of the criterion with some
parameter.
\begin{thm}\label{t1.6}
Let $u$ be a Leray-Hopf weak solution to the 3D Navier-Stokes
equations \eqref{a} with the initial value $u_{0}\in V\bigcap
L^{1}(\mathbb{R}^3)$. Suppose that, for some $i,j $ with $1\leq
i,j\leq 3$, $u$ satisfies one of the following conditions:
\par
(a) $i\neq j$,
\begin{eqnarray}\label{t15} \partial_iu_j\in L^{\beta_1}(0,T;
L^{\alpha_1}(\mathbb{R}^3))\ \mbox{with}\
\frac{3}{\alpha_1}+\frac{2}{\beta_1}=2-\epsilon,\end{eqnarray} where
\begin{eqnarray}\label{tt15} \displaystyle
\alpha_1\in[\frac{3}{2-\epsilon},\frac{3}{3-2\epsilon}],\
\mbox{and}\ \ 1\leq\epsilon<3/2.
\end{eqnarray}
\par
(b) $i=j$,
\begin{eqnarray}\label{ts15} \partial_ju_j\in L^{\beta_2}(0,T;
L^{\alpha_2}(\mathbb{R}^3))\  \mbox{with}\
\frac{3}{\alpha_2}+\frac{2}{\beta_2}=2-\epsilon,\end{eqnarray} where
\begin{eqnarray}\label{tts15} \displaystyle
\alpha_2\in[\frac{3}{2-\epsilon},\frac{6}{5-4\epsilon}], \ \mbox{
and}\  \ 1/2\leq\epsilon<5/4. \
\end{eqnarray} Then $u$ is regular.
\end{thm}
\begin{rem}
It is sufficient  to assume that $u_0\in V$ when we consider the
endpoint case of $\epsilon=1$ in part $(a)$ and $\epsilon=1/2 $ in
part $(b)$ respectively. In view of the result of  part $(b)$, we
can show the the line in Figure 1, which is continuous in
$\epsilon$. However, we see that this line is always under the
line``\textbf{(1)}" in Figure 1. From the proof of this Theorem, we
know that we choose an intermediate parameter $q$, and restrict $q$
to satisfy $1<q\leq2$ for convenience, which is the underlying
reason why the line is always below line``\textbf{(1)}" in Figure 1.
In fact, if we choose the intermediate parameter $q$ is larger than
2, we can get another better result such that the corresponding line
is always above line``\textbf{(1)}" in Figure 1, which is stated in
the following Theorem.
\end{rem}
\begin{thm}\label{t1.7}
Let $u$ be a Leray-Hopf weak solution to the 3D Navier-Stokes
equations \eqref{a} with the initial value $u_{0}\in V$. Suppose
that, for some $i,j $ with $1\leq i,j\leq 3$, $u$ satisfies one of
the following conditions:
\par
(a) $i\neq j$,
\begin{eqnarray}\label{jt15} \partial_iu_j\in L^{\beta_1}(0,T;
L^{\alpha_1}(\mathbb{R}^3))\ \mbox{with}\
\frac{3}{\alpha_1}+\frac{2}{\beta_1}=2-\epsilon,\end{eqnarray} where
\begin{eqnarray}\label{jtt15} \displaystyle
\alpha_1\in[\frac{3}{3-2\epsilon},\frac{3(11-2\epsilon)}{2\epsilon^2-26\epsilon+33}],
\ \mbox{and}\ \  1<\epsilon\leq21/16.
\end{eqnarray}
 \par $(b)$ $i=j$,
\begin{eqnarray}\label{tsw15} \partial_ju_j\in L^{\beta_2}(0,T;
L^{\alpha_2}(\mathbb{R}^3))\  \mbox{with}\
\frac{3}{\alpha_2}+\frac{2}{\beta_2}=2-\epsilon,\end{eqnarray} where
\begin{eqnarray}\label{ttsw15} \displaystyle
\alpha_2\in[\frac{6}{5-4\epsilon},\frac{18-2\epsilon}{(4\epsilon-3)(\epsilon-5)}],\
\mbox{and} \ \ 1/2<\epsilon\leq3/4.
\end{eqnarray} Then $u$ is regular.
\end{thm}
\begin{figure}
\centering
\includegraphics[width=0.8\textwidth,height=0.35\textheight]{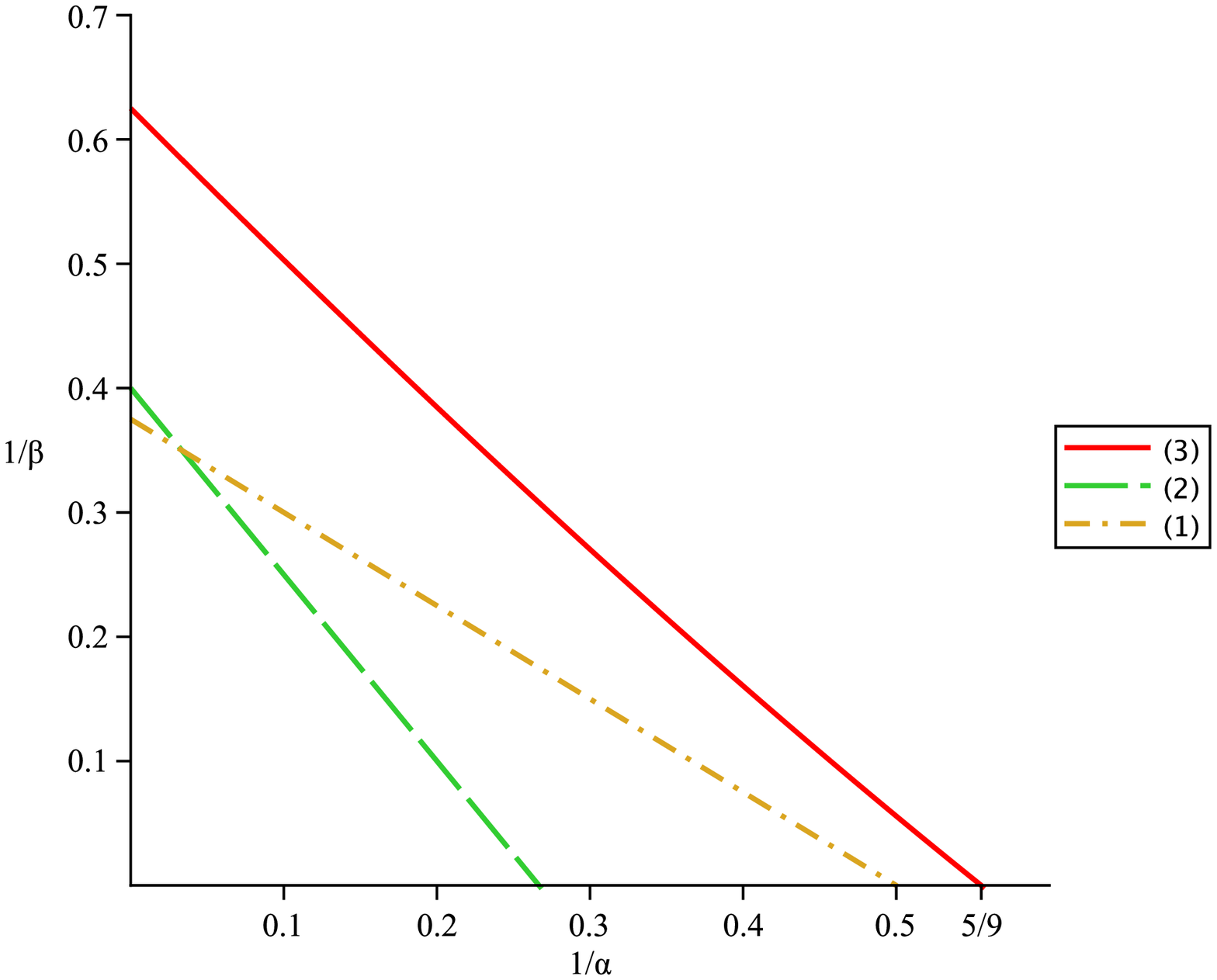}
\caption{Case of $j= k$} \label{fig:131}The line "\textbf{(1)}"
 is the result of  C.S. Cao, E.S.
Titi in \cite{[2]} ( see \eqref{d}). The line "\textbf{(3)}" is our
result, which mean \eqref{14}.  The result of Y. Zhou, M.
Pokorn$\acute{\mbox{y} }$ in \cite{[26]} (see \eqref{jbb}) is showed
by line "\textbf{(2)}".
\end{figure}
\par For  convenience, we recall the following version of
the three-dimensional  Sobolev and Ladyzhenskaya inequalities in the
whole space $\mathbb{R}^{3}$ (see, for example, \cite{[5]},
\cite{[8]}, \cite{[11]}). There exists a positive constant $C$ such
that
\begin{equation}\label{r}
\begin{array}{ll}
 \displaystyle
\|u\|_{r}&\displaystyle
\leq C \|u\|_{2}^{\frac{6-r}{2r}}\|\partial_{1}u\|_{2}^{\frac{r-2}{2r}}\|\partial_{2}u\|_{2}^{\frac{r-2}{2r}}\|\partial_{3}u\|_{2}^{\frac{r-2}{2r}}\\
&\leq C \|u\|_{2}^{\frac{6-r}{2r}}\|\nabla
u\|_{2}^{\frac{3(r-2)}{2r}},
\end{array}
\end{equation}
for every $u\in H^{1}(\mathbb{R}^{3})$ and every $r\in[2,6]$, where
$C$ is a constant depending only on $r$. Taking $\nabla$div on both
sides of \eqref{a} for smooth $(u; p)$, one can obtain
$$
-\Delta(\nabla
p)=\sum_{i,j}^{3}\partial_{i}\partial_{j}(\nabla(u_{i}u_{j})),
$$
therefore, the Calderon-Zygmund inequality in $\mathbb{R}^{3}$ (see
\cite{[30]})
\begin{equation}\label{s}
\|\nabla p\|_{q}\leq C\||\nabla u||u|\|_{q},
 1<q<\infty,
\end{equation}
holds, where $C$ is a positive constant depending only on $q$. And
there is another  estimates for the pressure
\begin{equation}\label{y33} \| p\|_{q}\leq C\|u\|_{2q}^2, \
 1<q<\infty,
\end{equation}

\section{{Proof of Main Results}}
In this section, under the assumptions of  Theorems
\ref{t1.3}-\ref{t1.4}, Theorem \ref{t1.9}, Theorem\ref{t1.5}
Theorem\ref{t1.6}, Theorem\ref{t1.7} in Section 1 respectively, we
prove our main results.  First of all, we note that, by the energy
inequality, for Leray-Hopf weak solutions, we have (see, for
example, \cite{[21]}, \cite{[22]} for detail)
\begin{equation}\label{e1}
\|u(\cdot,t)\|_{L^{2}}^{2}+2\nu\int_{0}^{t}\|\nabla
u(\cdot,s)\|_{L^{2}}^{2}ds\leq K_{1},\end{equation} for all $0<
t<T,$ where $K_{1}=\|u_{0}\|_{L^{2}}^{2}.$
\par  It is  well
known that there exists a unique strong solution  $u$ local in time
if $u_{0}\in V$. In addition, this strong solution $u\in
C([0,T^{*});V)\cap L^{2}(0,T^{*}; H^{2}(\mathbb{R}^{3}))$ is the
only weak solution with the initial datum $u_{0}$, where $(0,T^{*})$
is the maximal interval of existence of the unique strong solution.
If $T^{*}\geq T,$ then there is nothing to prove. If, on the other
hand, $T^{*}< T,$ then our strategy is to show that the $H^{1}$ norm
of this strong solution is bounded uniformly in time over the
interval $(0,T^{*})$, provided  additional conditions in Theorem
\ref{t1.3}-\ref{t1.4}, Theorem \ref{t1.9},
Theorem\ref{t1.5},Theorem\ref{t1.6}, Theorem\ref{t1.7}  in Section 1
are valid. As a result the interval $(0,T^{*})$ can not be a maximal
interval of existence, and consequently $T^{*}\geq T,$ which
concludes our proof. \par In order to prove the $H^{1}$ norm of the
strong solution $u$ is bounded on interval $(0,T^{*})$, combing with
the energy equality \eqref{e1}, it is sufficient to prove
\begin{equation} \label{e2}
\|\nabla
u\|_{L^2}^{2}+\displaystyle{\nu}\displaystyle\int_{0}^{t}\|\Delta
u\|_{L^2}^{2}d\tau\leq C,\ \forall \ t\in(0, T^{*})
\end{equation}
where  $C$  is a positive constant independent of $T^{*}.$ We recall
the following lemma (see \cite{[24]}), which is useful for our proof
of the Theorems.
\begin{lem}\label{l2.1}
Assume $v\in L^{\infty}_{t}L^{2}_{x}(\mathbb{R}^3\times I)$ and
$\nabla v\in L^{2}_{t}L^{2}_{x}(\mathbb{R}^3\times I)$, where $I$ is
an open interval. Then $v\in L^{s}_{t}L^{r}_{x}(\mathbb{R}^3\times
I)$ for all $r$ and $s$ such that
$$
\frac{2}{s}+\frac{3}{r}=\frac{3}{2}\ \mbox{with} \ 2\leq r\leq 6,
$$
and it holds $$ \|v\|_{ L^{s}_{t}L^{r}_{x}}\leq C
\|v\|_{L^{\infty}_{t}L^{2}_{x}}^{\frac{6-r}{2r}}\|\nabla
v\|_{L^{2}_{t}L^{2}_{x}}^{\frac{3(r-2)}{2r}}.
$$
\end{lem}
\noindent\textbf{Proof of Theorem \ref{t1.3}}  Taking the inner
product of the equation \eqref{a} with $-\Delta u$ in $L^{2}$, we
obtain
\begin{equation}\label{2227}\begin{array}{ll}
\displaystyle \frac{1}{2}\frac{d}{dt}\|\nabla
u\|_{2}^{2}&+\nu\|\Delta u\|_{2}^{2}\\&=\displaystyle
\sum_{i,j,k=1}^{3}\int_{\mathbb{R}^{3}}u_{i}\partial_{i}u_{j}\partial_{kk}u_{j}
dx\displaystyle \\ \displaystyle &=\displaystyle
\sum_{i,j,k=1}^{2}\int_{\mathbb{R}^{3}}u_{i}\partial_{i}u_{j}\partial_{kk}u_{j}dx
+ \sum_{i=1}^{3}
\sum_{j=1}^{2}\int_{\mathbb{R}^{3}}u_{i}\partial_{i}u_{j}\partial_{33}u_{j}dx\displaystyle \\
\displaystyle &\ \ \ \ +\displaystyle \sum_{j=1}^{2}
\sum_{k=1}^{2}\int_{\mathbb{R}^{3}}u_{3}\partial_{3}u_{j}\partial_{kk}u_{j}dx+\sum_{i,k=1}^{3}\int_{\mathbb{R}^{3}}u_{i}\partial_{i}u_{3}\partial_{kk}u_{3}dx\
\displaystyle \\
\displaystyle &=\displaystyle I_{1}(t)+I_{2}(t)+I_{3}(t)+I_{4}(t).
\end{array}\end{equation}
By integrating by parts a few times and using the incompressibility
condition, we have
$$
I_{1}(t)=\frac{1}{2}\sum_{i,j,k=1}^{2}\int_{\mathbb{R}^{3}}\partial_{i}u_{i}\partial_{k}u_{j}\partial_{k}u_{j}dx
-\sum_{i,j,k=1}^{2}\int_{\mathbb{R}^{3}}\partial_{k}u_{i}\partial_{i}u_{j}\partial_{k}u_{j}dx=I_{1}^{1}(t)+I_{1}^{2}(t).
$$
The terms $I_{1}^{1}(t),$ $ I_{1}^{2}(t), $ $I_{3}(t)$ and $I_{4}(t)
$ read as
$$
I_{1}^{1}(t)=-\frac{1}{2}\sum_{j,k=1}^{2}\int_{\mathbb{R}^{3}}\partial_{3}u_{3}\partial_{k}u_{j}\partial_{k}u_{j}dx,\hspace{4cm}
$$
$$\begin{array}{ll}\displaystyle
I_{1}^{2}(t)&=\displaystyle-\sum_{i,j,k=1}^{2}\int_{\mathbb{R}^{3}}\partial_{k}u_{i}\partial_{i}u_{j}\partial_{k}u_{j}dx\\
&
=\displaystyle\displaystyle\int_{\mathbb{R}^{3}}\partial_{2}u_{1}\partial_{1}u_{2}\partial_{2}u_{2}dx
+\displaystyle\int_{\mathbb{R}^{3}}\partial_{1}u_{2}\partial_{2}u_{1}\partial_{1}u_{1}dx
+\displaystyle\int_{\mathbb{R}^{3}}\partial_{1}u_{1}\partial_{1}u_{2}\partial_{1}u_{2}dx\\&
\ \
+\displaystyle\int_{\mathbb{R}^{3}}\partial_{1}u_{2}\partial_{2}u_{2}\partial_{1}u_{2}dx
+\displaystyle\int_{\mathbb{R}^{3}}\partial_{2}u_{1}\partial_{1}u_{1}\partial_{2}u_{1}dx
+\displaystyle\int_{\mathbb{R}^{3}}\partial_{1}u_{1}\partial_{1}u_{1}\partial_{1}u_{1}dx\\&\
\
+\displaystyle\int_{\mathbb{R}^{3}}\partial_{2}u_{2}\partial_{2}u_{2}\partial_{2}u_{2}dx
+\displaystyle\int_{\mathbb{R}^{3}}\partial_{2}u_{2}\partial_{2}u_{1}\partial_{2}u_{1}dx\\
&
=\displaystyle-\int_{\mathbb{R}^{3}}(\partial_{2}u_{1}\partial_{1}u_{2}\partial_{3}u_{3}
+\partial_{3}u_{3}\partial_{1}u_{2}\partial_{1}u_{2}+\partial_{2}u_{1}\partial_{3}u_{3}\partial_{2}u_{1})
dx\vspace{1mm}\\
&\ \
-\displaystyle\int_{\mathbb{R}^{3}}(\partial_{1}u_{1}\partial_{1}u_{1}\partial_{3}u_{3}
+\partial_{3}u_{3}\partial_{2}u_{2}\partial_{2}u_{2}-\partial_{1}u_{1}\partial_{3}u_{3}\partial_{2}u_{2})
dx,
\end{array}$$
$$\begin{array}{ll}\displaystyle
I_{3}(t)&=\displaystyle \sum_{j=1}^{2}
\sum_{k=1}^{2}\displaystyle\int_{\mathbb{R}^{3}}u_{3}\partial_{3}u_{j}\partial_{kk}u_{j}dx\\&=\displaystyle-\sum_{j=1}^{2}
\sum_{k=1}^{2}\int_{\mathbb{R}^{3}}\partial_{k}u_{3}\partial_{3}u_{j}\partial_{k}u_{j}dx-\sum_{j=1}^{2}
\sum_{k=1}^{2}\int_{\mathbb{R}^{3}}u_{3}\partial_{3k}u_{j}\partial_{k}u_{j}dx\\
&=\displaystyle-\sum_{j=1}^{2}
\sum_{k=1}^{2}\int_{\mathbb{R}^{3}}\partial_{k}u_{3}\partial_{3}u_{j}\partial_{k}u_{j}dx+\frac{1}{2}\displaystyle\sum_{j=1}^{2}
\sum_{k=1}^{2}\int_{\mathbb{R}^{3}}\partial_{3}u_{3}\partial_{k}u_{j}\partial_{k}u_{j}dx,
\end{array}$$
$$\begin{array}{ll}\displaystyle
I_{4}(t)&=\displaystyle
\sum_{i,k=1}^{3}\int_{\mathbb{R}^{3}}u_{i}\partial_{i}u_{3}\partial_{kk}u_{3}dx\\&=\displaystyle-
\sum_{i,k=1}^{3}\int_{\mathbb{R}^{3}}\partial_{k}u_{i}\partial_{i}u_{3}\partial_{k}u_{3}dx-
 \sum_{i,k=1}^{3}\int_{\mathbb{R}^{3}}u_{i}\partial_{ik}u_{3}\partial_{k}u_{3}dx\\
& =\displaystyle-
\sum_{i,k=1}^{3}\int_{\mathbb{R}^{3}}\partial_{k}u_{i}\partial_{i}u_{3}\partial_{k}u_{3}dx.
\end{array}$$
The above four  equalities imply  that
 $$
|I_{i}|\leq C \int_{\mathbb{R}^{3}}|u_{3}||\nabla u||\nabla^2 u|dx,\
i=1,3,4.
 $$
$\bullet\ \mbox{Case of} \
 \partial_{3}u_{2},\partial_{3}u_{3}$\\
As for $I_{2}$,  we have
$$\begin{array}{ll}\displaystyle
I_{2}&=-\displaystyle\int_{\mathbb{R}^{3}}\partial_{3}u_{1}\partial_{1}u_{1}\partial_{3}u_{1}
+\partial_{3}u_{2}\partial_{2}u_{1}\partial_{3}u_{1}+\partial_{3}u_{3}\partial_{3}u_{1}\partial_{3}u_{1}dx\vspace{1mm}\\
\displaystyle&\ \ \ \ \ \
-\displaystyle\int_{\mathbb{R}^{3}}\partial_{3}u_{1}\partial_{1}u_{2}\partial_{3}u_{2}+\partial_{3}u_{2}\partial_{2}u_{2}\partial_{3}u_{2}
+\partial_{3}u_{3}\partial_{3}u_{2}\partial_{3}u_{2}dx\vspace{1mm}\\
&\displaystyle=I_2^1+I_2^2+I_2^3+I_2^4+I_2^5+I_2^6.
\end{array}$$
It is obvious that \begin{equation} \label{e61}
\begin{array}{ll} \displaystyle
|I_2^j| \leq C \int_{\mathbb{R}^{3}}|u_{3}||\nabla u||\nabla^2
u|dx,\ j=3,6,
\end{array}\end{equation}
and
\begin{equation} \label{e62}
\begin{array}{ll} \displaystyle
|I_2^j| \leq C \int_{\mathbb{R}^{3}}|\partial_3u_{2}||\nabla
u|^2dx,\ j=2,4,5.
\end{array}\end{equation}
For the first term,  by using the incompressibility condition, and
integrating by parts a few times, we note that
\begin{equation} \label{e61}
\begin{array}{ll} \displaystyle
-\int_{\mathbb{R}^3}\partial_3u_1\partial_1u_1\partial_3u_1dx&=\displaystyle
2\int_{\mathbb{R}^3}\partial_{13}u_1\partial_3u_1u_1dx\vspace{1mm}\\
&=\displaystyle
-2\int_{\mathbb{R}^3}\partial_{23}u_2\partial_3u_1u_1dx-\displaystyle
2\int_{\mathbb{R}^3}\partial_{33}u_3\partial_3u_1u_1dx\vspace{1mm}\\
&=\displaystyle
2\int_{\mathbb{R}^3}\partial_{23}u_1\partial_3u_2u_1dx+\displaystyle
2\int_{\mathbb{R}^3}\partial_{3}u_2\partial_3u_1\partial_2u_1dx\vspace{1mm}\\
&\ \ \ +\displaystyle
2\int_{\mathbb{R}^3}\partial_{33}u_1\partial_3u_3u_1dx+\displaystyle
2\int_{\mathbb{R}^3}\partial_{3}u_3\partial_3u_1\partial_3u_1dx\vspace{1mm}\\
\end{array}\end{equation}
As before, $I_2^1$ has the estimate
\begin{equation} \label{e63}
\begin{array}{ll} \displaystyle
|I_2^1| &\leq \displaystyle C
\int_{\mathbb{R}^{3}}|\partial_3u_{2}||\nabla
u|^2dx+C\int_{\mathbb{R}^{3}}|u_{3}||\nabla u||\nabla^2
u|dx\vspace{1mm}\\&\ \ \displaystyle
+C\int_{\mathbb{R}^{3}}|\partial_3u_{2}||\nabla^2
u||u|dx+C\int_{\mathbb{R}^{3}}|\partial_3u_{3}||\nabla^2 u||u|dx.
\end{array}\end{equation}
Therefore, by above inequalities, we see that
\begin{equation} \label{12h}\begin{array}{ll}\displaystyle \frac{1}{2}\frac{d}{dt}\|\nabla
u\|_{2}^{2}+\nu\|\Delta u\|_{2}^{2}&\leq \displaystyle C
\int_{\mathbb{R}^{3}}|\partial_3u_{2}||\nabla
u|^2dx+C\int_{\mathbb{R}^{3}}|u_{3}||\nabla u||\nabla^2
u|dx\vspace{1mm}\\&\ \ \displaystyle
+C\int_{\mathbb{R}^{3}}|\partial_3u_{2}||\nabla^2
u||u|dx+C\int_{\mathbb{R}^{3}}|\partial_3u_{3}||\nabla^2 u||u|dx.
\end{array}\end{equation}
\par Let $T^{\prime}\in(0,T)$ with $T^\prime\leq T^{\ast}$ be
arbitrary. We shall prove that $T^\prime$ is not a blow-up point. By
decreasing of $\alpha_i, i=1,2,$ if necessary, we may assume
\begin{eqnarray}\label{e49} u_3\in L^{\beta_1}(0,T;
L^{\alpha_1}(\mathbb{R}^3))\ \mbox{with}\
\frac{3}{\alpha_1}+\frac{2}{\beta_1}=1, \alpha_1\in(3,\infty], 
\end{eqnarray}
and
\begin{eqnarray}\label{e50} \partial_3u_i\in L^{\beta_2}(0,T;
L^{\alpha_2}(\mathbb{R}^3))\ \mbox{with}\
\frac{3}{\alpha_2}+\frac{2}{\beta_2}=2, \alpha_2\in[2,3], 
 i=2,3.
\end{eqnarray}
Choose a $t_1\in(0,T^\prime)$ such that
\begin{eqnarray}\label{e51} \|\partial_3u_i\|_{ L^{\beta_2}(t_1,T^{\prime};
L^{\alpha_2}(\mathbb{R}^3))}\leq\epsilon,\ i=2,3,
\end{eqnarray}
and
\begin{eqnarray}\label{e52} \|u_3\|_{ L^{\beta_1}(t_1,T^{\prime};
L^{\alpha_1}(\mathbb{R}^3))}\leq\epsilon,
\end{eqnarray}
where $\epsilon$ is to be determined. Let $t_2\in(t_1,T^\prime)$ be
arbitrary. For any $r,s,$ we abbreviate
\begin{eqnarray}\label{e53} \|\cdot\|_{L^{s}_{t}L^{r}_{x}}=\|\cdot\|_{ L^{s}((t_1,t_2);
L^{r}(\mathbb{R}^3))}.
\end{eqnarray}We denote
\begin{equation}\label{e47}
L=\displaystyle\|\nabla u\|_{L_{t}^{\infty}L_{x}^{2}}+\nu\|\Delta
u\|_{L_{t}^{2}L_{x}^{2}}.
\end{equation}
Let $t\in(t_1, t_2]$ be arbitrary.  Integrating \eqref{12h} on
$(t_1, t)$, we get
\begin{equation}\label{e62}\begin{array}{ll}
\|\nabla u\|_{L^2}^{2}&+\displaystyle
2\nu\displaystyle\int_{t_1}^{t}\|\Delta u\|_{L^2}^{2}d\tau\\
&\leq\displaystyle C
\int_{t_1}^{t}\int_{\mathbb{R}^{3}}|\partial_3u_{2}||\nabla
u|^2dxd\tau+C \int_{t_1}^{t}\int_{\mathbb{R}^{3}}|u_{3}||\nabla
u||\nabla^2 u|dxd\tau\vspace{1mm}\\&\ \ \displaystyle +C
\int_{t_1}^{t}\int_{\mathbb{R}^{3}}|\partial_3u_{2}||\nabla^2
u||u|dxd\tau+C
\int_{t_1}^{t}\int_{\mathbb{R}^{3}}|\partial_3u_{3}||\nabla^2
u||u|dxd\tau\vspace{1mm}\\
&\ \  \displaystyle+\|\nabla
u(t_1)\|_{L^2}^{2}\vspace{1mm}\\&=\displaystyle L_1+L_2+L_3+L_4+L_5.
\end{array}
\end{equation}
We estimate $L_{i},$ $i=1,2,3,4$ one by one. Firstly, for $L_1$ we
have
\begin{equation} \label{e63}\begin{array}{ll}
L_1&\leq\displaystyle
C\|\partial_3u_2\|_{L_{t}^{\beta_2}L_{x}^{\alpha_2}}\|\nabla u\|_{L_{t}^{s_1}L_{x}^{r_1}}^2\vspace{1mm}\\
\end{array}\end{equation}
where $s_1$ and $r_1$ satisfy
\begin{equation} \label{e64}
\frac{1}{\alpha_2}+\frac{2}{r_1}=1,\ \
\frac{1}{\beta_2}+\frac{2}{s_1}=1,
\end{equation}
by \eqref{e50} and \eqref{e64}, we have $3\leq r_1\leq4$ and
$2/s_1+3/r_1=3/2.$ By Lemma 2.1, we have
\begin{equation} \label{e65}
L_1\leq\displaystyle C\epsilon L^2.\end{equation} As for $L_2$, we
have
\begin{equation}\label{e56}\begin{array}{ll}
L_2&\leq\displaystyle C\|u_3\|_{L_{t}^{\beta_1}L_{x}^{\alpha_1}}
\|\Delta u\|_{L_{t}^{2}L_{x}^{2}}\|\nabla u\|_{L_{t}^{s_2}L_{x}^{r_2}}\vspace{2mm}\\
&\leq\displaystyle C\epsilon L^2,
\end{array}\end{equation}
where $s_2$ and $r_2$ satisfy
\begin{equation}\label{e57}
\frac{1}{\alpha_1}+\frac{1}{r_2}=\frac{1}{2},\ \
\frac{1}{\beta_1}+\frac{1}{s_2}=\frac{1}{2},
\end{equation}
by \eqref{e49} and \eqref{e57}, we have $2\leq r_2<6$ and
$2/s_2+3/r_2=3/2.$ thus, $s_2$ and $r_2$ satisfy Lemma 2.1. The
estimate of $L_3$ is as follows
\begin{equation} \label{e68}\begin{array}{ll}
L_3&\leq\displaystyle
C\|\partial_3u_2\|_{L_{t}^{\beta_2}L_{x}^{\alpha_2}} \|\Delta
u\|_{L_{t}^{2}L_{x}^{2}}
\|u\|_{L_{t}^{s_3}L_{x}^{r_3}}\vspace{1mm}\\
\end{array}\end{equation}
where $s_3$ and $r_3$ satisfy
\begin{equation}\label{e58}
\frac{1}{\alpha_2}+\frac{1}{r_3}=\frac{1}{2},\ \
\frac{1}{\beta_2}+\frac{1}{s_3}=\frac{1}{2},
\end{equation}
by \eqref{e50} and \eqref{e58}, we have
\begin{equation}\label{e59}
\frac{2}{s_3}+\frac{3}{r_3}=\frac{1}{2},\  r_3\geq6.\end{equation}By
the Gagliardo-Nirenberg inequality
\begin{equation} \label{e588}
\|v\|_{L^r}\leq\displaystyle\|\nabla v\|_{L^{\frac{3r}{3+r}}},\
\frac{3}{2}\leq r<\infty.\end{equation} Therefore, we have
\begin{equation} \label{e589}
\|u\|_{L_{t}^{s_3}L_{x}^{r_3}}\leq\displaystyle\|\nabla
u\|_{L_{t}^{s_3}L_{x}^{r_4}},\end{equation} where
$r_4=3r_3/(3+r_3)$, and we have $2/s_3+3/r_4=3/2$ with $2\leq
r_4\leq 6.$ Combining \eqref{e50}, \eqref{e68} and \eqref{e589}, as
well as Lemma 2.1, one has
\begin{equation} \label{e65}
L_3\leq\displaystyle C\epsilon L^2.\end{equation} The term $L_4$ is
estimated the same way and we get the same result. Finally, we
obtain
$$
L^2\leq C\epsilon L^2+ \|\nabla u(t_1)\|_{L^2}^{2}.
$$
If $t_1$ is sufficiently close to $T^\prime$ and $\epsilon$ is
sufficiently small, we can absorb the first term into  the left hand
side, and then we obtain that $L$ is bounded with a bound
independent of $t_2\in(t_1,T^\prime)$. Finally, we get $\|\nabla
u\|_{L_{t}^{\infty}L_{x}^{2}(\mathbb{R}^3\times(t_1, T^\prime))}\leq
C$. Therefore, the solution cannot blow up at $T^\prime$. We
complete
the proof of the case of $\partial_{3}u_{2},\partial_{3}u_{3}$.\\
 $\bullet\ \mbox{Case of} \
 \partial_{3}u_{1},\partial_{3}u_{3}$\par This case is similar to
 the first case. The main difference is that we have another estimate for
$L_2$.  Also by using the incompressibility condition, and
integrating by parts a few times, we note that $I_2^5$ becomes
\begin{equation} \label{e61}
\begin{array}{ll} \displaystyle
-\int_{\mathbb{R}^3}\partial_3u_2\partial_3u_2\partial_2u_2dx&=\displaystyle
2\int_{\mathbb{R}^3}\partial_{23}u_2\partial_3u_2u_2dx\vspace{1mm}\\
&=\displaystyle
-2\int_{\mathbb{R}^3}\partial_{13}u_1\partial_3u_2u_2dx-\displaystyle
2\int_{\mathbb{R}^3}\partial_{33}u_3\partial_3u_2u_2dx\vspace{1mm}\\
&=\displaystyle
2\int_{\mathbb{R}^3}\partial_{13}u_2\partial_3u_1u_2dx+\displaystyle
2\int_{\mathbb{R}^3}\partial_{3}u_2\partial_3u_1\partial_1u_2dx\vspace{1mm}\\
&\ \ \ +\displaystyle
2\int_{\mathbb{R}^3}\partial_{33}u_2\partial_3u_3u_2dx-\displaystyle
2\int_{\mathbb{R}^3}\partial_{3}u_3\partial_3u_2\partial_3u_2dx,
\end{array}\end{equation}
and then $I_2^5$ has the estimate
\begin{equation} \label{e63}
\begin{array}{ll} \displaystyle
|I_2^5| &\leq \displaystyle C
\int_{\mathbb{R}^{3}}|\partial_3u_{1}||\nabla
u|^2dx+C\int_{\mathbb{R}^{3}}|u_{3}||\nabla u||\nabla^2
u|dx\vspace{1mm}\\&\ \ \displaystyle
+C\int_{\mathbb{R}^{3}}|\partial_3u_{1}||\nabla^2
u||u|dx+C\int_{\mathbb{R}^{3}}|\partial_3u_{3}||\nabla^2 u||u|dx.
\end{array}\end{equation}
Therefore, we finally get
\begin{equation} \label{192}\begin{array}{ll}\displaystyle \frac{1}{2}\frac{d}{dt}\|\nabla
u\|_{2}^{2}+\nu\|\Delta u\|_{2}^{2}&\leq \displaystyle C
\int_{\mathbb{R}^{3}}|\partial_3u_{1}||\nabla
u|^2dx+C\int_{\mathbb{R}^{3}}|u_{3}||\nabla u||\nabla^2
u|dx\vspace{1mm}\\&\ \ \displaystyle
+C\int_{\mathbb{R}^{3}}|\partial_3u_{1}||\nabla^2
u||u|dx+C\int_{\mathbb{R}^{3}}|\partial_3u_{3}||\nabla^2 u||u|dx.
\end{array}\end{equation}
By using \eqref{192}, we give the same method as before to get the
desired
result. The proof is completed. \\
\textbf{Proof of Theorem
\ref{t1.4}} From the condition of this Theorem, we split
the proof into two parts. \\
$\bullet\  i\neq j$
\par Firstly, we consider the case of $q=2,$ and then we see that the range of $\alpha_2$ is
$\alpha_2>2$. For convenience of writing, we set
\begin{equation}\label{e18}
r=\frac{3\alpha_2-2}{\alpha_2}.
\end{equation}
It is easy to check that $r>2$ when $\alpha_2>2$. Without loss of
generality, in the proof, we will assume that $i=1, j=2,$ the other
cases can be discussed in the same way (for details see Remark 2.2
below). We begin with \eqref{2227}, and the same process to the
proof of Theorem \ref{t1.3}, we firstly have
 $$
|I_{i}|\leq C \int_{\mathbb{R}^{3}}|u_{3}||\nabla u||\nabla^2 u|dx,\
i=1,3,4.
 $$
As for $I_{2}$, also by the incompressibility condition,  we have
$$\begin{array}{ll}\displaystyle
I_{2}&=-\displaystyle\int_{\mathbb{R}^{3}}\partial_{3}u_{1}\partial_{1}u_{1}\partial_{3}u_{1}
+\partial_{3}u_{2}\partial_{2}u_{1}\partial_{3}u_{1}+\partial_{3}u_{3}\partial_{3}u_{1}\partial_{3}u_{1}dx\vspace{1mm}\\
\displaystyle&\ \ \ \ \ \
-\displaystyle\int_{\mathbb{R}^{3}}\partial_{3}u_{1}\partial_{1}u_{2}\partial_{3}u_{2}+\partial_{3}u_{2}\partial_{2}u_{2}\partial_{3}u_{2}
+\partial_{3}u_{3}\partial_{3}u_{2}\partial_{3}u_{2}dx\vspace{1mm}\\
&\displaystyle=-\displaystyle\int_{\mathbb{R}^{3}}\partial_{3}u_{1}\left(-\partial_{2}u_{2}-\partial_{3}u_{3}\right)\partial_{3}u_{1}
+\partial_{3}u_{2}\partial_{2}u_{1}\partial_{3}u_{1}+\partial_{3}u_{3}\partial_{3}u_{1}\partial_{3}u_{1}dx\vspace{1mm}\\
\displaystyle&\ \ \ \ \ \
-\displaystyle\int_{\mathbb{R}^{3}}\partial_{3}u_{1}\partial_{1}u_{2}\partial_{3}u_{2}+\partial_{3}u_{2}\partial_{2}u_{2}\partial_{3}u_{2}
+\partial_{3}u_{3}\partial_{3}u_{2}\partial_{3}u_{2}dx\vspace{1mm}\\
&\displaystyle\leq\displaystyle\int_{\mathbb{R}^{3}}|u_{2}||\nabla
u||\nabla^2 u|dx +\int_{\mathbb{R}^{3}}|u_{3}||\nabla u||\nabla^2
u|dx.
\end{array}$$
Therefore, we get
\begin{equation}\label{3337}\begin{array}{ll}
\displaystyle &\displaystyle\frac{1}{2}\frac{d}{dt}\|
\nabla u\|_{L^2}^{2}+\nu\|\Delta u\|_{L^2}^{2}\vspace{1mm}\\
&\ \ \ \ \ \ \ \ \leq C\displaystyle
\int_{\mathbb{R}^{3}}|u_3||\nabla u||\nabla^2 u|dx+
C\int_{\mathbb{R}^{3}}|u_2||\nabla u||\nabla^2u|dx,\vspace{1mm}\\
&\ \ \ \ \ \ \ \ =K_1(t)+K_2(t).
\end{array}\end{equation}
Next, we estimate $K_1(t)$ and $K_2(t)$. Firstly, we pay attention
to $K_2(t)$, applying H$\ddot{\mbox{o}}$lder's inequality several
times, we obtain
\begin{equation} \label{e19}\begin{array}{ll}
K_{2}(t) \displaystyle &\leq\displaystyle
 C\int_{\mathbb{R}^{2}}\max_{x_{1}}|u_{2}|(\int_{\mathbb{R}}
 |\nabla u|^{2}dx_{1})^{\frac{1}{2}}(\int_{\mathbb{R}}|\nabla^2 u|^{2}dx_{1})^{\frac{1}{2}} dx_{2}dx_{3}\
 \displaystyle \vspace{2mm}\\
 \displaystyle &\ \ \ \ \leq\displaystyle
 C[\int_{\mathbb{R}^{2}}(\max_{x_{1}}|u_{2}|)^{r}dx_{2}dx_{3}]^{\frac{1}{r}}
 [\int_{\mathbb{R}^{2}}(\int_{\mathbb{R}}|\nabla
 u|^{2}dx_{1})^{\frac{r}{r-2}}dx_{2}dx_{3}]^{\frac{r-2}{2r}}\\
 \displaystyle& \ \ \ \ \ \ \times
 \displaystyle[\int_{\mathbb{R}^{3}}|\nabla^2
 u|^{2}dx_{1}dx_{2}dx_{3}]^{\frac{1}{2}}\\
 \displaystyle \vspace{2mm}
 \displaystyle &\ \ \ \ \leq\displaystyle
 C[\int_{\mathbb{R}^{3}}|u_{2}|^{r-1}|\partial_{1}u_{2}|dx_{1}dx_{2}dx_{3}]^{\frac{1}{r}}\|\Delta u\|_{2}
  \displaystyle \\
 \displaystyle& \ \ \ \ \ \ \times \displaystyle[\int_{\mathbb{R}}(\int_{\mathbb{R}^{2}}|\nabla
 u|^{\frac{2r}{r-2}}dx_{2}dx_{3})^{\frac{r-2}{r}}dx_{1}]^{\frac{1}{2}}
 \displaystyle \vspace{2mm}\\
\displaystyle &\ \ \ \ \ \leq \displaystyle C\displaystyle\|
u_{2}\|_{L^2}^{\frac{r-1}{r}}\|
\partial_{1}u_{2}\|_{L^\frac{2}{3-r}}^{\frac{1}{r}}\|\nabla u\|_{L^2}^{\frac{r-2}{r}}
 \|\partial_{2}\nabla u\|_{L^2}^{\frac{1}{r}} \|\partial_{3}\nabla u\|_{L^2}^{\frac{1}{r}}
\|\Delta u\|_{2}\vspace{2mm}\\
\displaystyle &\ \ \ \ \ \leq \displaystyle C\displaystyle\|
u_{2}\|_{L^2}^{\frac{r-1}{r}}\|
\partial_{1}u_{2}\|_{L^\frac{2}{3-r}}^{\frac{1}{r}}\|\nabla u\|_{L^2}^{\frac{r-2}{r}}
\|\Delta u\|_{L^2}^{\frac{r+2}{r}}.
\end{array}\end{equation}
In above inequality, from \eqref{e18}, we note that
$\frac{2}{3-r}=\alpha_2. $ Therefore, applying  Young's inequality,
\eqref{e19} immediately implies
\begin{equation} \label{e20}\begin{array}{ll}
K_{2}(t)\displaystyle &\leq \displaystyle C\displaystyle\|
u_{2}\|_{L^2}^{\frac{2(\alpha_2-1)}{3\alpha_2-2}}\|
\partial_{1}u_{2}\|_{L^{\alpha_2}}^{\frac{\alpha_2}{3\alpha_2-2}}\|\nabla u\|_{L^2}^{\frac{\alpha_2-2}{3\alpha_2-2}}
\|\Delta u\|_{L^2}^{\frac{5\alpha_2-2}{3\alpha_2-2}}\vspace{2mm}\\
&\leq \displaystyle C\displaystyle\|
u_2\|_{L^2}^{\frac{4(\alpha_2-1)}{\alpha_2-2}}\|
\partial_{1}u_{2}\|_{L^{\alpha_2}}^{\frac{2\alpha_2}{\alpha_2-2}}\|\nabla
u\|_{L^2}^{2}+\frac{\nu}{4} \|\Delta u\|_{L^2}^{2}.
\end{array}\end{equation}
As for $K_1(t)$,  applying H$\ddot{\mbox{o}}$lder's and Young's
inequalities, we have
\begin{equation} \label{28}\begin{array}{ll}
K_{1}(t)& \displaystyle =\displaystyle C
\int_{\mathbb{R}^{3}}|u_3||\nabla
u||\nabla^2 u|dx\vspace{1mm}\\
&\displaystyle\leq C\|u_3\|_{L^{\alpha_1}}\|\nabla u\|_{L^q}\|\Delta u\|_{L^2}\vspace{1mm}\\
&\displaystyle\leq C \|u_3\|_{L^{\alpha_1}}\|\nabla
u\|_{{L^2}}^{\frac{6-q}{2q}}\|\Delta
u\|_{{L^2}}^{\frac{5q-6}{2q}}\vspace{1mm}\\
&\displaystyle\leq C \|u_3\|_{L^{\alpha_1}}^{\frac{4q}{6-q}}\|\nabla
u\|_{{L^2}}^{2}+\frac{\nu}{4}\|\Delta u\|_{{L^2}}^{2},
\end{array}\end{equation}
where $\alpha_1$ and $q$ satisfy
\begin{equation}\label{e15}
\frac{1}{q}+\frac{1}{\alpha_1}=\frac{1}{2} \ \mbox{with}\ 2\leq q<6,
\alpha_1>3.
\end{equation}
From \eqref{28} and \eqref{e20}, one has
\begin{equation} \label{e21}
\begin{array}{ll} \displaystyle \frac{1}{2}\frac{d}{dt}\|
\nabla u\|_{L^2}^{2}+\nu\|\Delta u\|_{L^2}^{2} &\displaystyle\leq
C\| u_{2}\|_{L^2}^{\frac{4(\alpha_2-1)}{\alpha_2-2}}\|
\partial_{1}u_{2}\|_{L^\alpha_2}^{\frac{2\alpha_2}{\alpha_2-2}}\|
\nabla u\|_{L^2}^{2} \vspace{2mm}\\
&\displaystyle \ \ \ \ \ +C\displaystyle\|
u_3\|_{L^\alpha_1}^{\frac{2\alpha_1}{\alpha_1-3}}\|\nabla
u\|_{L^2}^{2}+\frac{\nu}{2} \|\Delta u\|_{L^2}^{2}.\vspace{2mm}
\end{array}\end{equation}
Absorbing the last term in right hand of \eqref{e21} and integrating
the inequality on time, using the energy inequality, we obtain
\begin{equation} \label{e22}\begin{array}{ll}
\displaystyle\| \nabla u\|_{L^2}^{2}+\nu\int_{0}^{t}\|\nabla
\Delta u\|_{L^2}^{2}d\tau\displaystyle&\leq \displaystyle
C\int_{0}^{t}\| u_{2}\|_{L^2}^{\frac{4(\alpha_2-1)}{\alpha_2-2}}\|
\partial_{1}u_{2}\|_{L^{\alpha_2}}^{\frac{2\alpha_2}{\alpha_2-2}}\|\nabla u\|_{L^2}^2d\tau\vspace{1mm}\\
&\displaystyle \ \ \ +C\int_{0}^{t}\|u_3\|_{L^{\alpha_1}}^{\frac{2\alpha_1}{\alpha_1-3}}\|\nabla u\|_{L^2}^2d\tau
+\|\nabla u(0)\|_{L^2}^2\vspace{1mm}\\
&\leq \displaystyle C\int_{0}^{t}\|
\partial_{1}u_{2}\|_{L^{\alpha_2}}^{\frac{2\alpha_2}{\alpha_2-2}}\|\nabla u\|_{L^2}^2d\tau\vspace{1mm}\\
&\displaystyle
\ \ \ +C\int_{0}^{t}\|u_3\|_{L^{\alpha_1}}^{\frac{2\alpha_1}{\alpha_1-3}}
\|\nabla u\|_{L^2}^2d\tau+\|\nabla u(0)\|_{L^2}^2.\vspace{1mm}\
\end{array}\end{equation}
By  using  Gronwall's inequality, we obtain
\begin{equation}\label{e23}\begin{array}{ll} \displaystyle &\displaystyle\|
\nabla u\|_{L^2}^{2}+\nu\int_{0}^{t}\|\Delta u\|_{L^2}^{2}d\tau\vspace{1mm}\\
&\ \ \ \ \ \ \ \ \displaystyle\leq
\left(\|\nabla u(0)\|_{L^2}^2\right)
\exp\left(C\int_{0}^{t}\|\partial_1u_2\|_{L^{\alpha_2}}^{\frac{2\alpha_2}{\alpha_2-2}}
d\tau\right)\exp\left(C\int_{0}^{t}\|u_3\|_{L^{\alpha_1}}^{\frac{2\alpha_1}{\alpha_1-3}}
d\tau\right)
\end{array}\end{equation}
By condition \eqref{15}$-$\eqref{17},  \eqref{e23} follows that the
$ H^1$ norm of the strong solution $u$ is bounded on the maximal
interval of existence $(0, T^{*})$.\par
 Now, we pay attention to the case of  $1<q<2$. Next, we give an estimate on
$u_2$. We use $|u_{2}|^{q-1}\mbox{sgn}(u_2)$ with $1<q\leq3/2$ as
test function in the equation \eqref{a} for $u_{2}.$ By using
Gagliardo-Nirenberg and H$\ddot{\mbox{o}}$lder's inequalities
H$\ddot{\mbox{o}}$lder's inequalities and \eqref{s}, we have
\begin{equation}\label{td9}
\begin{array}{ll}
 \displaystyle \frac{1}{q}\frac{d}{dt}\|u_{2}\|_{L^q}^{q}&
 +C(q)\nu\|\nabla|u_{2}|^{\frac{q}{2}}\|_{L^2}^{2}\\
 &=\displaystyle
 -\int_{\mathbb{R}^{3}}\partial_{2}p|u_{2}|^{q-1}\mbox{sgn}(u_2)dx\displaystyle  \vspace{2mm}\\
 \displaystyle &\leq\displaystyle
  C\|\nabla p\|_{L^q}\|u_{2}\|_{L^q}^{q-1}\ \vspace{2mm}\\
 \displaystyle &\leq\displaystyle   C\||\nabla u||u|\|_{L^q}\|u_{2}\|_{L^q}^{q-1}\vspace{2mm}\\
 \displaystyle &\leq\displaystyle C\|\nabla u\|_{L^2}\| u\|_{L^{\frac{2q}{2-q}}}
\|u_{2}\|_{L^q}^{q-1}
 \vspace{2mm}\\
\displaystyle &\leq\displaystyle C\|\nabla u\|_{L^2}\|
u\|_{L^2}^{\frac{3-2q}{q}}\|\nabla u\|_{L^2}^{\frac{3q-3}{q}}
 \|u_{2}\|_{L^q}^{q-1}\vspace{2mm}\\
\displaystyle &=\displaystyle C\| u\|_{L^2}^{\frac{3-2q}{q}}\|\nabla
u\|_{L^2}^{\frac{4q-3}{q}}
 \|u_{2}\|_{L^q}^{q-1},
\end{array}
\end{equation}
where we note that $1<q\leq3/2$ means $2<\frac{2q}{2-q}\leq 6$, and
\eqref{td9} immediately  implies that
$$
\frac{d}{dt}\|u_{2}\|_{L^q}\leq C\|
u\|_{L^2}^{\frac{3-2q}{q}}\|\nabla u\|_{L^2}^{\frac{4q-3}{q}}.
$$
After integrating on time, and note that $u_0\in V\bigcap
L^{q}(\mathbb{R}^3)$, by energy and H$\ddot{\mbox{o}}$lder's
inequalities one has
\begin{equation}\label{ta10}
\begin{array}{ll}
 \displaystyle
\|u_{2}\|_{L^q}&\leq\|u_{2}(0)\|_{L^q}+\displaystyle C\int_{0}^{t}
\| u\|_{L^2}^{\frac{3-2q}{q}}\|\nabla
u\|_{L^2}^{\frac{4q-3}{q}}d\tau\\
&\leq\|u_{2}(0)\|_{L^q}+\displaystyle C\int_{0}^{t} \|\nabla
u\|_{L^2}^{\frac{4q-3}{q}}d\tau\\
&\leq\|u_{2}(0)\|_{L^q}+\displaystyle C\left(\int_{0}^{t} \|\nabla
u\|_{L^2}^{2}d\tau\right)^{\frac{4q-3}{2q}}T^{\frac{3-2q}{2q}}\\
&\leq\|u_{3}(0)\|_{L^q}+\displaystyle C(T).
\end{array}\end{equation}
Therefore, we get
\begin{equation}\label{td13}
u_2\in L^{\infty}(0, T; L^{q}(\mathbb{R}^{3}))\ \mbox{with}\
1<q\leq\frac{3}{2}.
\end{equation}
On the other hand, by energy inequality we know that
\begin{equation}\label{td14}
u_2\in L^{\infty}(0, T; L^{2}(\mathbb{R}^{3})).
\end{equation}
Finally, by interpolation, we have
\begin{equation}\label{td15}
u_2\in L^{\infty}(0, T; L^{q}(\mathbb{R}^{3}))\ \mbox{with}\ 1<q<2.
\end{equation}
For every $\alpha_2>\frac{q}{q-1}$, we set
$$
r=\frac{(q+1)\alpha_2-q}{\alpha_2},
$$
then we have $r>2$. Similar to \eqref{e19}, also by
H$\ddot{\mbox{o}}$lder's and Young's inequalities, we obtain another
estimate
\begin{equation} \label{te41}\begin{array}{ll}
K_{2}(t) \displaystyle & \leq \displaystyle C\displaystyle\|
u_{2}\|_{L^q}^{\frac{r-1}{r}}\|
\partial_{1}u_{2}\|_{L^\frac{q}{q-r+1}}^{\frac{1}{r}}\|\nabla u\|_{L^2}^{\frac{r-2}{r}}
\|\Delta u\|_{L^2}^{\frac{r+2}{r}}\vspace{2mm}\\
&\leq \displaystyle C\displaystyle\|
u_{2}\|_{L^q}^{\frac{2(r-1)}{r-2}}\|
\partial_{1}u_{2}\|_{L^{\alpha_2}}^{\frac{2}{r-2}}
\|\nabla u\|_{L^2}^{2}+\frac{\nu}{4}\|\Delta u\|_{L^2}^{2}.\vspace{2mm}\\\
\end{array}\end{equation}
Applying \eqref{td15} and integrating above inequality, one has
\begin{equation} \label{te42}\begin{array}{ll}
\displaystyle\displaystyle\int_0^tK_{2}(\tau)d\tau &\leq
\displaystyle C\displaystyle\int_0^t\displaystyle\|
u_{2}\|_{L^q}^{\frac{2(r-1)}{r-2}}\|
\partial_{1}u_{2}\|_{L^{\alpha_2}}^{\frac{2}{r-2}}
\|\nabla u\|_{L^2}^{2}d\tau+\frac{\nu}{4}\int_0^t\|\Delta u\|_{L^2}^{2}d\tau\vspace{2mm}\\
&\leq \displaystyle C\displaystyle\int_0^t\displaystyle\|
\partial_{1}u_{2}\|_{L^{\alpha_2}}^{\frac{2}{r-2}}
\|\nabla u\|_{L^2}^{2}d\tau+\frac{\nu}{4}\int_0^t\| \Delta u\|_{L^2}^{2}d\tau.\vspace{2mm}\\
\end{array}\end{equation}
Integrating \eqref{3337} on time, and absorbing the last term in
\eqref{te42} and \eqref{28} respectively, it follows that
\begin{equation} \label{te43}\begin{array}{ll}
\displaystyle &\displaystyle\|
\nabla u\|_{L^2}^{2}+\nu\int_{0}^{t}\|\Delta u\|_{L^2}^{2}d\tau\vspace{1mm}\\
&\ \ \ \ \ \ \ \ \displaystyle\leq
C\int_{0}^{t}\|\partial_1u_2\|_{L^{\alpha_2}}^{\frac{2}{r-2}}\|\nabla
u\|_{L^2}^2d\tau
+C\int_{0}^{t}\|u_3\|_{L^{\alpha_1}}^{\frac{2\alpha_1}{2\alpha_1-3}}\|\nabla u\|_{L^2}^2d\tau+\|\nabla u(0)\|_{L^2}^2.\\
\end{array}\end{equation}
By using Gronwall's inequality and condition
\eqref{15}$-$\eqref{17}, we also obtain $ H^1$ norm of the strong
solution $u$ is bounded on the maximal interval of existence $(0,
T^{*})$ when $1<q<2$. Thus we prove $(a)$.\\
$\bullet \ i=j$
\par Without loss of generality, here, we assume $i=j=2$. Similar to the
proof of the part $(a)$, we estimate the second term $K_2(t)$ of
\eqref{3337}. Firstly, for every
\begin{equation}\label{e32}\frac{9}{5}<\alpha_3\leq\infty,\end{equation} we choose
\begin{equation}\label{e31}\left\{\begin{array}{l}
\displaystyle \frac{1}{\gamma}=\frac{3+\sqrt{24\alpha_3^2-24\alpha_3+9}}{12\alpha_3},\\
\displaystyle \mu=\frac{\alpha_3\gamma}{2\alpha_3-\gamma}.\vspace{2mm}\\
\end{array}
\right.
\end{equation}
From \eqref{e31}$_1$ and $\eqref{e32}$, we have
$\frac{1}{\gamma}>\frac{1}{2\alpha_3}$, which means that
$\gamma<2\alpha_3$, and hance $\mu>0$ is well defined in
\eqref{e31}$_2$. Also from \eqref{e31}$_1$,  we see that $\gamma$ is
an increasing function  with the variable $\alpha_3\in(9/5,\infty]$,
and from $\eqref{e32}$
 we get
\begin{equation}\label{e33}2<\frac{9}{4}<\gamma\leq\sqrt{6},\end{equation}
 and moreover, \eqref{e31}
follows
\begin{equation}\label{e34} \frac{1}{\mu}+\frac{1}{\alpha_3}+\frac{\gamma-2}{\gamma}=1.\end{equation}
Besides, again by \eqref{e31}, we see that
\begin{equation}\label{e29}
\mu=\frac{6\alpha_3}{-3+\sqrt{24\alpha_3^2-24\alpha_3+9}},
\end{equation}
by the monotonicity,  we obtain
\begin{equation}\label{e35}
1<\frac{\sqrt 6}{2}\leq\mu<3.
\end{equation}
We choose
\begin{equation}\label{e36}
\beta_3=\frac{2\mu}{3-\mu},
\end{equation}
then we have
$$\frac{1}{\mu}=\frac{2}{3}(\frac{1}{\beta_3}+\frac{1}{2}),$$
by \eqref{e31}$_1$, \eqref{e34} and \eqref{e36}, we can compute that
\begin{equation}\label{e37}
\frac{2}{\beta_3}+1=\frac{3}{\mu}=
\frac{6}{\gamma}-\frac{3}{\alpha_3}=\frac{-3+\sqrt{24\alpha_3^2-24\alpha_3+9}}{2\alpha_3}\Rightarrow
\frac{3}{2\alpha_3}+\frac{2}{\beta_3}=f(\alpha_3).\end{equation}
 \par Now, we use $|u_{2}|^{\gamma-2}u_{2}$ with
$\gamma>2$ as test function in the equation \eqref{a} for $u_{2}.$
By using of Gagliardo-Nirenberg and H$\ddot{\mbox{o}}$lder's
inequalities, we have
\begin{equation}\label{e24}
\begin{array}{ll}
 \displaystyle \frac{1}{\gamma}\frac{d}{dt}\|u_{2}\|_{L^\gamma}^{\gamma}&
 +C(\gamma)\nu\|\nabla|u_{2}|^{\frac{\gamma}{2}}\|_{L^2}^{2}\\
 &=\displaystyle
 -\int_{\mathbb{R}^{3}}\partial_{2}p|u_{2}|^{\gamma-2}u_{2}dx\displaystyle  \vspace{2mm}\\
 &\leq C\displaystyle\int_{\mathbb{R}^{3}}|p||u_{2}|^{\gamma-2}|\partial_{2}u_{2}|dx\displaystyle  \vspace{2mm}\\
 \displaystyle &\leq\displaystyle
  C\|p\|_{L^\mu}\|u_{2}\|_{L^\gamma}^{\gamma-2}\|\partial_{2}u_{2}\|_{{L^{\alpha_3}}}\ \vspace{2mm}\\
 \displaystyle &\leq\displaystyle   C \|u
 \|_{L^{2\mu}}^2\|u_{2}\|_{L^\gamma}^{\gamma-2}\|\partial_2u_{2}\|_{L^{\alpha_3}}\displaystyle
 \ (\mbox{by }\eqref{y33})\vspace{2mm}\\
 \displaystyle &\leq\displaystyle C\| u\|_{L^2}^{\frac{3-\mu}{\mu}}\|\nabla
 u\|_{L^2}^{\frac{3(\mu-1)}{\mu}}
 \|u_{2}\|_{L^\gamma}^{\gamma-2}\|\partial_{2}u_{2}\|_{L^{\alpha_3}},
\end{array}
\end{equation}
where the $\gamma, \mu$ and $\alpha_3$ satisfy \eqref{e34}. The above
inequality immediately implies that
\begin{equation}\label{e26}
\begin{array}{ll}
\displaystyle
\frac{1}{2}\frac{d}{dt}\|u_{2}\|_{L^\gamma}^{2}\displaystyle
\leq\displaystyle C\| u\|_{L^2}^{\frac{3-\mu}{\mu}}\|\nabla
 u\|_{L^2}^{\frac{3(\mu-1)}{\mu}}\|\partial_{2}u_{2}\|_{L^{\alpha_3}}.
\end{array}
\end{equation}
In view of  \eqref{e35}, we have $\frac{3(\mu-1)}{\mu}<2$, applying
Young's inequality, we have
\begin{equation}\label{e27}
\begin{array}{ll}
\displaystyle
\frac{1}{2}\frac{d}{dt}\|u_{2}\|_{L^\gamma}^{2}\displaystyle
\leq\displaystyle C\|\nabla
 u\|_{L^2}^{2}+\| u\|_{L^2}^{2}\|\partial_{2}u_{2}\|_{L^{\alpha_3}}^{\frac{2\mu}{3-\mu}}.
\end{array}
\end{equation}
Integrating \eqref{e27} on time, and by energy inequality
\eqref{e1}, we obtain
\begin{equation}\label{e38}
\begin{array}{ll}
\displaystyle \|u_{2}\|_{L^\gamma}^{2}\displaystyle &\leq\displaystyle
\|u_{2}(0)\|_{L^\gamma}^{2}+C+C\int_{0}^{t}\|\partial_{2}u_{2}\|_{L^{\alpha_3}}^{\frac{2\mu}{3-\mu}}d\tau\vspace{1mm}\\
 &=\displaystyle
\|u_{2}(0)\|_{L^\gamma}^{2}+C+C\int_{0}^{t}\|\partial_{2}u_{2}\|_{L^{\alpha_3}}^{\beta_3}d\tau.
\end{array}
\end{equation}
In view of $\eqref{e33}$, we have $\|u_{2}(0)\|_{\gamma}<C$ for some
$C>0$, therefore, by the condition \eqref{9}, we get
\begin{equation}\label{e39}
u_2\in L^{\infty}(0, T; L^{\gamma}(\mathbb{R}^{3})).
\end{equation}
For the mentioned parameter $\gamma$ in \eqref{e31}, we set
\begin{equation}\label{e40}
r=\frac{(\gamma+1)\alpha_3-\gamma}{\alpha_3},
\end{equation}
then
$$
r=\frac{12\alpha_3-9+\sqrt{24\alpha_3^2-24\alpha_3+9}}{3+\sqrt{24\alpha_3^2-24\alpha_3+9}},
$$
since $9/5<\alpha_3\leq\infty$, also by monotonicity, it is easy to
see $2<r\leq\sqrt 6+1$. Similar to \eqref{e19}, also by
H$\ddot{\mbox{o}}$lder's and Young's inequalities, we obtain another
estimate
\begin{equation} \label{e41}\begin{array}{ll}
K_{2}(t) \displaystyle & \leq \displaystyle C\displaystyle\|
u_{2}\|_{L^\gamma}^{\frac{r-1}{r}}\|
\partial_{2}u_{2}\|_{L^\frac{\gamma}{\gamma-r+1}}^{\frac{1}{r}}\|\nabla u\|_{L^2}^{\frac{r-2}{r}}
\|\Delta u\|_{L^2}^{\frac{r+2}{r}}\vspace{2mm}\\
&\leq \displaystyle C\displaystyle\|
u_{2}\|_{L^\gamma}^{\frac{2(r-1)}{r-2}}\|
\partial_{2}u_{2}\|_{L^{\alpha_3}}^{\frac{2}{r-2}}
\|\nabla u\|_{L^2}^{2}+\frac{\nu}{4}\|\Delta u\|_{L^2}^{2}.
\end{array}\end{equation}
Applying \eqref{e39} and integrating above inequality, one has
\begin{equation} \label{e41}\begin{array}{ll}
\displaystyle\displaystyle\int_0^tK_{2}(\tau)d\tau &\leq
\displaystyle C\displaystyle\int_0^t\|
u_{2}\|_{L^\gamma}^{\frac{2(r-1)}{r-2}}\|
\partial_{2}u_{2}\|_{L^{\alpha_3}}^{\frac{2}{r-2}}
\|\nabla u\|_{L^2}^{2}d\tau+\frac{\nu}{4}\int_0^t\|\Delta u\|_{L^2}^{2}d\tau\vspace{2mm}\\
&\leq \displaystyle C\displaystyle\int_0^t\|
\partial_{2}u_{2}\|_{L^{\alpha_3}}^{\frac{2}{r-2}}
\|\nabla u\|_{L^2}^{2}d\tau+\frac{\nu}{4}\int_0^t\| \Delta
u\|_{L^2}^{2}d\tau.
\end{array}\end{equation}
Integrating \eqref{3337} on time, and absorbing the last term in
\eqref{e41} and \eqref{28} respectively, it follows that
\begin{equation} \label{e43}\begin{array}{ll}
\displaystyle &\displaystyle\|
\nabla u\|_{L^2}^{2}+\nu\int_{0}^{t}\|\Delta u\|_{L^2}^{2}d\tau\vspace{1mm}\\
&\ \ \ \ \ \ \ \ \displaystyle\leq
C\int_{0}^{t}\|\partial_2u_2\|_{L^{\alpha_3}}^{\frac{2}{r-2}}\|\nabla
u\|_{L^2}^2d\tau
+C\int_{0}^{t}\|u_3\|_{L^{\alpha_1}}^{\frac{2\alpha_1}{2\alpha_1-3}}\|\nabla
u\|_{L^2}^2d\tau+\|\nabla u(0)\|_{L^2}^2.
\end{array}\end{equation}
We claim $\frac{2}{r-2}=\beta_3$. In fact, from \eqref{e31}, we have
$$
\gamma=\frac{-3+\sqrt{24\alpha_3^2-24\alpha_3+9}}{2(\alpha_3-1)},
$$
and then by the definition \eqref{e40}, one has
\begin{equation}\label{e44}
\frac{2}{\frac{2}{r-2}}+1=r-1=\gamma-\frac{\gamma}{\alpha_3}=\frac{-3+\sqrt{24\alpha_3^2-24\alpha_3+9}}{2\alpha_3}.
\end{equation}
Comparing \eqref{e44} with \eqref{e37}, we prove the claim.
Therefore, we can apply  Gronwall's inequality to \eqref{e43}, and
by condition \eqref{15}, \eqref{8} and \eqref{9} to get
that the $ H^1$ norm of the strong solution $u$ is bounded on the
maximal interval of existence $(0, T^{*})$. The proof of this
Theorem is completed.
\begin{rem}
In the process of the proof, if we want to prove the case of $i=3$
when $j=2$.  For the first part, the inequality \eqref{e19} may be
replaced by
\begin{equation} \label{}\begin{array}{ll}
K_{2}(t)\leq \displaystyle C\displaystyle\|
u_{2}\|_{L^2}^{\frac{r-1}{r}}\|
\partial_{3}u_{2}\|_{L^\frac{2}{3-r}}^{\frac{1}{r}}\|\nabla u\|_{L^2}^{\frac{r-2}{r}}
\|\Delta u\|_{L^2}^{\frac{r+2}{r}},
\end{array}\end{equation}
and \eqref{e20} becomes
\begin{equation} \label{}\begin{array}{ll}
K_{2}(t)\displaystyle &\leq \displaystyle C\displaystyle\|
u_{2}\|_{L^2}^{\frac{4(\alpha_2-1)}{\alpha_2-2}}\|
\partial_{3}u_{2}\|_{L^{\alpha_2}}^{\frac{2\alpha_2}{\alpha_2-2}}\|
\nabla u\|_{L^2}^{2}+\frac{\nu}{4}
\|\Delta u\|_{L^2}^{2}.\vspace{2mm}\\
\end{array}\end{equation}
If we want to prove the case of $j=1$, we shall  give an alternative
proof of the term $I_2$, also by the incompressibility condition,
one has
$$\begin{array}{ll}\displaystyle
I_{2}&=-\displaystyle\int_{\mathbb{R}^{3}}\partial_{3}u_{1}\partial_{1}u_{1}\partial_{3}u_{1}
+\partial_{3}u_{2}\partial_{2}u_{1}\partial_{3}u_{1}+\partial_{3}u_{3}\partial_{3}u_{1}\partial_{3}u_{1}dx\vspace{1mm}\\
\displaystyle&\ \ \ \ \ \
-\displaystyle\int_{\mathbb{R}^{3}}\partial_{3}u_{1}\partial_{1}u_{2}\partial_{3}u_{2}+\partial_{3}u_{2}\partial_{2}u_{2}\partial_{3}u_{2}
+\partial_{3}u_{3}\partial_{3}u_{2}\partial_{3}u_{2}dx\vspace{1mm}\\
&\displaystyle=-\displaystyle\int_{\mathbb{R}^{3}}\partial_{3}u_{1}\partial_{1}u_{1}\partial_{3}u_{1}
+\partial_{3}u_{2}\partial_{2}u_{1}\partial_{3}u_{1}+\partial_{3}u_{3}\partial_{3}u_{1}\partial_{3}u_{1}dx\vspace{1mm}\\
\displaystyle&\ \ \ \ \ \
-\displaystyle\int_{\mathbb{R}^{3}}\partial_{3}u_{1}\partial_{1}u_{2}\partial_{3}u_{2}+\partial_{3}u_{2}\left(-\partial_{1}u_{1}-\partial_{3}u_{3}\right)\partial_{3}u_{2}
+\partial_{3}u_{3}\partial_{3}u_{2}\partial_{3}u_{2}dx\vspace{1mm}\\
&\displaystyle\leq\displaystyle\int_{\mathbb{R}^{3}}|u_{3}||\nabla u||\nabla^2 u|dx+\int_{\mathbb{R}^{3}}|u_{1}||\nabla u||\nabla^2 u|dx.\vspace{1mm}\\
\end{array}$$
and then we obtain
\begin{equation}\label{4447}\begin{array}{ll}
\displaystyle &\displaystyle\frac{1}{2}\frac{d}{dt}\|
\nabla u\|_{L^2}^{2}+\nu\|\Delta u\|_{L^2}^{2}\vspace{1mm}\\
&\ \ \ \ \ \ \ \ \leq C\displaystyle
\int_{\mathbb{R}^{3}}|u_3||\nabla u||\nabla^2 u|dx+
C\int_{\mathbb{R}^{3}}|u_1||\nabla u||\nabla^2u|dx,\vspace{1mm}\\
&\ \ \ \ \ \ \ \ =K_{1}(t)+K_{2}(t),
\end{array}\end{equation} by which one can prove the case of $i\neq j=1$   and $i=j=1.$  The term $K_{1}(t)$ is the same to \eqref{28}. As
for  the second term $K_{2}(t)$, we shall give the inequality
\eqref{e19} as the following for $i\neq j=1$,
\begin{equation} \label{}\begin{array}{ll}
K_{2}(t)\leq \displaystyle C\displaystyle\|
u_{1}\|_{L^2}^{\frac{r-1}{r}}\|
\partial_{i}u_{1}\|_{L^\frac{2}{3-r}}^{\frac{1}{r}}\|\nabla u\|_{L^2}^{\frac{r-2}{r}}
\|\Delta u\|_{L^2}^{\frac{r+2}{r}}, \ i=2 \ \mbox{or}\ 3,
\end{array}\end{equation}
and then get the corresponding form of \eqref{e20}. As for $i=j=1$,
we will use $|u_{1}|^{\gamma-2}u_{1},$ as test function in the
equation for $u_{1},$  and we can get the similar results to
\eqref{e24}. By the same process to prove the case of $i=j=1.$
\end{rem}
\noindent\textbf{Proof of Theorem \ref{t1.9}}   Give the same
process as in the  Theorem \ref{t1.3}, we have
\begin{equation} \label{12}\begin{array}{ll}\displaystyle \frac{1}{2}\frac{d}{dt}\|\nabla
u\|_{L^2}^{2}+\nu\|\Delta u\|_{L^2}^{2}&\leq \displaystyle C
\int_{\mathbb{R}^{3}}|\partial_3u_{2}||\nabla
u|^2dx+C\int_{\mathbb{R}^{3}}|u_{3}||\nabla u||\nabla^2
u|dx\vspace{1mm}\\&\ \ \displaystyle
+C\int_{\mathbb{R}^{3}}|\partial_3u_{2}||\Delta
u||u|dx+C\int_{\mathbb{R}^{3}}|\partial_3u_{3}||\Delta u||u|dx
\vspace{1mm}\\&\displaystyle =L_1+L_2+L_3+L_4.
\end{array}\end{equation}
Here, we only prove the case of $\partial_3u_2, \partial_3u_3$. For
the case of $\partial_3u_1, \partial_3u_3$, we will begin with
\eqref{192}, and from which we can give the similar proof.\par We
estimate $L_{i},$ $i=1,2,3,4$ one by one. Firstly, for $L_1$ we have
\begin{equation} \label{p8}\begin{array}{ll}
L_1&\leq\displaystyle C\|\partial_3u_2\|_{L_{}^{\alpha_1}}\|\nabla
u\|_{L_{}^{r_1}}^2,
\vspace{1mm}\\
\end{array}\end{equation}
where  $r_1$ satisfies
\begin{equation} \label{p9}
\frac{1}{\alpha_1}+\frac{2}{r_1}=1\Rightarrow
r_1=\frac{2\alpha_1}{\alpha_1-1},
\end{equation}
since $2\leq\alpha_1\leq3$, we have $3\leq r_1\leq4$, by
Gagliardo-Nirenberg and Young's inequalities, on has
\begin{equation} \label{p10}\begin{array}{ll}
L_1&\leq\displaystyle C\|\partial_3u_2\|_{L_{}^{\alpha_1}}\|\nabla
u\|_{L_{}^{2}}^{\frac{6-r_1}{r_1}}\|\Delta
u\|_{L_{}^{2}}^{\frac{3r_1-6}{r_1}}\vspace{1mm}\\
&\leq\displaystyle
C\|\partial_3u_2\|_{L_{}^{\alpha_1}}^{\frac{2r_1}{6-r_1}}\|\nabla u\|_{L_{}^{2}}^{2}+\frac{\nu}{4}\|\Delta u\|_{L_{}^{2}}^{2}\vspace{1mm}\\
&=\displaystyle
C\|\partial_3u_2\|_{L_{}^{\alpha_1}}^{\frac{2\alpha_1}{2\alpha_1-3}}\|\nabla u\|_{L_{}^{2}}^{2}+\frac{\nu}{4}\|\Delta u\|_{L_{}^{2}}^{2}.\vspace{1mm}\\
\end{array}\end{equation}
As for $L_3$, let   $r_2$ satisfy
\begin{equation}\label{p11}
\frac{1}{\alpha_1}+\frac{1}{r_2}=\frac{1}{2},\ \
\end{equation} then we have
\begin{equation} \label{p12}\begin{array}{ll}
L_3&\leq\displaystyle C\|\partial_3u_2\|_{L_{}^{\alpha_1}} \|\Delta
u\|_{L_{}^{2}} \|u\|_{L_{}^{r_2}}.\vspace{1mm}\\
\end{array}\end{equation}
 From the fact that $2\leq\alpha_1\leq 3$, we have $r_2\geq 6.$
By the Gagliardo-Nirenberg inequality,
\begin{equation} \label{p13}
\|v\|_{L^r}\leq\displaystyle\|\nabla v\|_{L^{\frac{3r}{3+r}}},\
\frac{3}{2}\leq r<\infty.\end{equation} Therefore, we have
\begin{equation} \label{p14}
\|u\|_{L_{}^{r_2}}\leq\displaystyle\|\nabla
u\|_{L_{}^{r_3}},\end{equation} where $r_3=3r_2/(3+r_2)$ with $2\leq
r_3\leq 3.$ Thus, we have

\begin{equation}\label{p15}
\frac{1}{\alpha_1}+\frac{1}{r_3}=\frac{5}{6}\Rightarrow
r_3=\frac{6\alpha_1}{5\alpha_1-6},\
\end{equation}
and applying Young's inequality, \eqref{p12} becomes
\begin{equation} \label{p16}\begin{array}{ll}
L_3 &\leq\displaystyle C\|\partial_3u_2\|_{L_{}^{\alpha_1}} \|\Delta
u\|_{L_{}^{2}} \|\nabla u\|_{L_{}^{r_3}}.\vspace{1mm}\\
&\leq\displaystyle C\|\partial_3u_2\|_{L_{}^{\alpha_1}} \|\nabla
u\|_{L_{}^{2}}^{\frac{6-r_3}{2r_3}} \|\Delta u\|_{L_{}^{2}}^{\frac{5r_3-6}{2r_3}}.\vspace{1mm}\\
&\leq\displaystyle
C\|\partial_3u_2\|_{L_{}^{\alpha_1}}^{\frac{4r_3}{6-r_3}} \|\nabla
u\|_{L_{}^{2}}^{2}+\frac{\nu}{8} \|\Delta u\|_{L_{}^{2}}^{2}.\vspace{1mm}\\
&=\displaystyle
C\|\partial_3u_2\|_{L_{}^{\alpha_1}}^{\frac{2\alpha_1}{2\alpha_1-3}}
\|\nabla
u\|_{L_{}^{2}}^{2}+\frac{\nu}{8} \|\Delta u\|_{L_{}^{2}}^{2}.\vspace{1mm}\\
\end{array}\end{equation}
The term $L_4$ is estimated in the same way and we get
\begin{equation} \label{p17}\begin{array}{ll}
L_4 &\leq\displaystyle
C\|\partial_3u_3\|_{L_{}^{\alpha_2}}^{\frac{2\alpha_2}{2\alpha_2-3}}
\|\nabla
u\|_{L_{}^{2}}^{2}+\frac{\nu}{8} \|\Delta u\|_{L_{}^{2}}^{2}.\vspace{1mm}\\
\end{array}\end{equation}
Therefore, integrating on time and absorbing the last term in
\eqref{p10}, \eqref{p16} and \eqref{p17}, we get
\begin{equation} \label{p18}\begin{array}{ll}
&\|\nabla
u\|_{L^2}^{2}+\displaystyle\frac{3\nu}{2}\displaystyle\int_{0}^{t}\|\Delta
u\|_{L^2}^{2}d\tau\\
&\ \ \ \ \ \  \leq\|\nabla u(0)\|_{L^2}^{2}+\displaystyle
C\int_{0}^{t}\|
\partial_{3}u_{2}\|_{L^{\alpha_1}}^{\frac{2\alpha_1}{2\alpha_1-3}}\|\nabla
u\|_{L^2}^{2}d\tau\\
&\hspace{1.3cm} +C\displaystyle\int_{0}^{t}\|
\partial_{3}u_{3}\|_{L^{\alpha_2}}^{\frac{2\alpha_2}{2\alpha_2-3}}\|\nabla
u\|_{L^2}^{2}d\tau+C\displaystyle\int_{0}^{t} L_{2}(\tau)d\tau.
\end{array}\end{equation}

As for the  estimation of  $L_2$, we give the same proof as the case
of $i=j$ in Theorem \ref{t1.4}, in which $u_2$ is replaced by $u_3$,
$\partial_2u_2$ is replaced by $\partial_3u_3$ and $\alpha_3,
\beta_3$ is replaced by $\alpha_2,\beta_2$. Finally, we get
\begin{equation} \label{p37}\begin{array}{ll}
\displaystyle\displaystyle\int_0^tL_{2}(\tau)d\tau &\leq
\displaystyle C\displaystyle\int_0^t\|
u_{3}\|_{L^{\gamma}}^{\frac{2(r-1)}{r-2}}\|
\partial_{3}u_{3}\|_{L^{\alpha_2}}^{\frac{2}{r-2}}
\|\nabla u\|_{L^2}^{2}d\tau+\frac{\nu}{2}\int_0^t\|\Delta u\|_{L^2}^{2}d\tau\vspace{2mm}\\
&\leq \displaystyle C\displaystyle\int_0^t\|
\partial_{3}u_{3}\|_{L^{\alpha_2}}^{\frac{2}{r-2}}
\|\nabla u\|_{L^2}^{2}d\tau+\frac{\nu}{2}\int_0^t\| \Delta u\|_{L^2}^{2}d\tau,\vspace{2mm}\\
\end{array}\end{equation}
where $\alpha_2\geq2$ and
\begin{equation}\label{p21}\left\{\begin{array}{l}
\displaystyle \frac{1}{\gamma}=\frac{3+\sqrt{24\alpha_2^2-24\alpha_2+9}}{12\alpha_2},\\
\displaystyle \mu=\frac{\alpha_2\gamma}{2\alpha_2-\gamma},\vspace{2mm}\\
\displaystyle r=\frac{(\gamma+1)\alpha_2-\gamma}{\alpha_2}.
\end{array}
\right.
\end{equation}
Inserting \eqref{p37} into \eqref{p18} and  absorbing the last term
in \eqref{p37}, note the boundedness of $\| u_{3}\|_{\gamma}$, we
have
\begin{equation} \label{p38}\begin{array}{ll}
&\|\nabla
u\|_{L^2}^{2}+\displaystyle{\nu}\displaystyle\int_{0}^{t}\|\Delta
u\|_{L^2}^{2}d\tau\\
&\ \ \ \ \ \  \leq\|\nabla u(0)\|_{L^2}^{2}+\displaystyle
C\int_{0}^{t}\|
\partial_{3}u_{2}\|_{L^{\alpha_1}}^{\frac{2\alpha_1}{2\alpha_1-3}}\|\nabla
u\|_{L^2}^{2}d\tau\\
&\hspace{1.3cm} +C\displaystyle\int_{0}^{t}\|
\partial_{3}u_{3}\|_{L^{\alpha_2}}^{\frac{2\alpha_2}{2\alpha_2-3}}\|\nabla
u\|_{L^2}^{2}d\tau+C\displaystyle \displaystyle\int_0^t\|
\partial_{3}u_{3}\|_{L^{\alpha_2}}^{\frac{2}{r-2}}
\|\nabla u\|_{L^2}^{2}d\tau.
\end{array}\end{equation}
We also can check that $\frac{2}{r-2}=\beta_2$ (see \eqref{e44}),
where $\beta_2=\frac{2\mu}{3-\mu}$. Denote that
$$
\frac{1}{\beta^\prime}=\frac{2\alpha_2-3}{2\alpha_2}.
$$
Also from \eqref{e44}, $
r=\frac{-3+\sqrt{24\alpha_2^2-24\alpha_2+9}}{2\alpha_2}+1, $ we see
that
\begin{equation} \label{p42}\begin{array}{ll}
\displaystyle\displaystyle \frac{1}{\beta^\prime}- \frac{1}{\beta_2}
&\displaystyle
\displaystyle=\frac{2\alpha_2-3}{2\alpha_2}-\frac{r-2}{2}\vspace{2mm}\\&
\displaystyle \displaystyle=\frac{(4-r)\alpha_2-3}{2\alpha_2}\vspace{2mm}\\
&=\displaystyle \displaystyle\frac{\frac{6\alpha_2+3-\sqrt{24\alpha_2^2-24\alpha_2+9}}{2}-3}{2\alpha_2}\vspace{2mm}\\
&=\displaystyle \displaystyle\frac{6\alpha_2-3-\sqrt{24\alpha_2^2-24\alpha_2+9}}{4\alpha_2}\vspace{2mm}\\
&=\displaystyle \displaystyle\frac{12\alpha_2(\alpha_2-1)}
{4\alpha_2\left((6\alpha_2-3)+\sqrt{24\alpha_2^2-24\alpha_2+9}\right)}>0,\ \forall \ \alpha_2\geq2.\vspace{2mm}\\
\end{array}\end{equation}
From above inequality, we know  $\beta^\prime<\beta_2$ and hance
$\beta^\prime/\beta_2<1$. Therefore, applying
H$\ddot{\mbox{o}}$lder's inequality, and by energy inequality, one
has
\begin{equation} \label{p43}\begin{array}{ll}
\displaystyle\displaystyle\int_{0}^{t}\|
\partial_{3}u_{3}\|_{L^{\alpha_2}}^{\frac{2\alpha_2}{2\alpha_2-3}}\|\nabla
u\|_{L^2}^{2}d\tau &=\displaystyle \displaystyle\int_{0}^{t}\|
\partial_{3}u_{3}\|_{L^{\alpha_2}}^{\beta^\prime}\|\nabla
u\|_{L^2}^{2}d\tau\vspace{2mm}\\
&=\displaystyle \displaystyle\int_{0}^{t}\|
\partial_{3}u_{3}\|_{L^{\alpha_2}}^{\beta^\prime}\|\nabla
u\|_{L^2}^{\frac{2\beta^\prime}{\beta_2}}d\tau\|\nabla
u\|_{L^2}^{\frac{2(\beta_2-\beta^\prime)}{\beta_2}}d\tau\vspace{2mm}\\
&\leq \displaystyle C\displaystyle\left(\int_0^t\|
\partial_{3}u_{3}\|_{L^{\alpha_2}}^{\beta_2}
\|\nabla
u\|_{2}^{2}d\tau\right)^{\frac{\beta^\prime}{\beta_2}}\left(\int_0^t\|
\nabla u\|_{L^2}^{2}d\tau\right)
^{\frac{\beta_2-\beta^\prime}{\beta_2}}\vspace{2mm}\\
&\leq \displaystyle C\displaystyle\left(\int_0^t\|
\partial_{3}u_{3}\|_{L^{\alpha_2}}^{\beta_2}
\|\nabla
u\|_{L^2}^{2}d\tau\right)^{\frac{\beta^\prime}{\beta_2}}\vspace{2mm}\\
&\leq \displaystyle C\displaystyle\int_0^t\|
\partial_{3}u_{3}\|_{L^{\alpha_2}}^{\frac{2}{r-2}}
\|\nabla u\|_{L^2}^{2}d\tau+C.
\end{array}\end{equation}
Finally, we get
\begin{equation} \label{p44}\begin{array}{ll}
\|\nabla
u\|_{L^2}^{2}&+\displaystyle{\nu}\displaystyle\int_{0}^{t}\|\Delta
u\|_{L^2}^{2}d\tau\\
&\ \ \ \ \ \  \leq\|\nabla u(0)\|_{L^2}^{2}+\displaystyle
C\int_{0}^{t}\|
\partial_{3}u_{2}\|_{L^{\alpha_1}}^{\frac{2\alpha_1}{2\alpha_1-3}}\|\nabla
u\|_{2}^{2}d\tau\\
&\hspace{1.3cm} +C\displaystyle \displaystyle\int_0^t\|
\partial_{3}u_{3}\|_{L^{\alpha_2}}^{\frac{2}{r-2}}
\|\nabla u\|_{L^2}^{2}d\tau+C.
\end{array}\end{equation}
By using Gronwall's inequality to \eqref{p44}, and by condition
\eqref{o6}, \eqref{o8} and \eqref{o9}, we get that the $ H^1$ norm
of the strong solution $u$ is bounded on the maximal interval of
existence $(0, T^{*})$. The proof of the case of
$\partial_{3}u_{2},\partial_{3}u_{3}$ is completed.\\
\noindent\textbf{Proof of Theorem \ref{t1.5}}  Without loss of
generality, we
 assume $j=3, k=3$. For every
\begin{equation}\label{d1} \alpha\in(\frac{9}{5},\infty],\end{equation}
we take
\begin{equation}\label{d2}\left\{\begin{array}{l}
\displaystyle \mu=\frac{24\alpha}{\alpha-12+\sqrt{144-264\alpha+289\alpha^2}},\vspace{2mm}\\
\displaystyle q_2=\frac{2\alpha\mu}{\alpha+\mu},\vspace{2mm}\\
\displaystyle r=\frac{2\mu\alpha-\mu+\alpha}{\alpha+\mu}.
\end{array}
\right.
\end{equation}
From \eqref{d2}, we see that $\mu$ is a decreasing function of
$\alpha$, and by \eqref{d1},  we have
\begin{equation}\label{d5}\frac{4}{3}\leq\mu<3,\end{equation}
and
\begin{equation}\label{d6} \frac{1}{\mu}+\frac{1}{\alpha}+\frac{q_2-2}{q_2}=1.\end{equation}
On the other hand, from \eqref{d2}, we have
$$
q_2=\frac{48\alpha}{\alpha+12+\sqrt{144-264\alpha+289\alpha^2}},
$$
and $q_2$ is an increasing function with  $\alpha\in(9/5,\infty]$,
which follows
$$9/4< q_2\leq 8/3.$$ We choose
\begin{equation}\label{d7}
\beta=\frac{2\mu}{3-\mu},
\end{equation}
then from \eqref{d2}, \eqref{d6} and \eqref{d7}, we can compute that
\begin{equation}\label{e377}
\frac{2}{\beta}+1=\frac{3}{\mu}=\frac{\alpha+\sqrt{289\alpha^2-264\alpha+144}}{8\alpha}-\frac{3}{2\alpha}\Rightarrow
\frac{3}{2\alpha}+\frac{2}{\beta}=g(\alpha).\end{equation} We denote
\begin{equation}\label{d8}
V_{1}(t)=\int_{0}^{t}\|
\partial_{3}u_{3}\|_{L^\alpha}^{\beta}\|\nabla
u\|_{L^2}^{2}d\tau=\int_{0}^{t}\|
\partial_{3}u_{3}\|_{L^\alpha}^{\frac{2\mu}{3-\mu}}\|\nabla
u\|_{L^2}^{2}d\tau.
\end{equation} Next, we give a estimate of
$u_3$. We use $|u_{3}|^{q_2-2}u_{3}$  as a test function  in the
equation \eqref{a} for $u_{3}.$ Similar to \eqref{e24}, we have
\begin{equation}\label{d12}
\begin{array}{ll}
\displaystyle \|u_{3}\|_{L^{q_2}}^{2}\displaystyle \leq\displaystyle
\|u_{3}(0)\|_{L^{q_2}}^{2}+C+C\int_{0}^{t}\|\partial_{3}u_{3}\|_{L^\alpha}^{\frac{2\mu}{3-\mu}}d\tau.
\end{array}
\end{equation}
By the condition \eqref{14}, $\eqref{d2}$ and \eqref{d5}, we have
$q_2<6$. Note that $\|u_{3}(0)\|_{q_2}<C$ for some $C>0$, we get
\begin{equation}\label{d13}
u_3\in L^{\infty}(0, T; L^{q_2}(\mathbb{R}^{3})).
\end{equation}
By \eqref{d2}, we have
$$r=\frac{(q_{2}+1){\alpha}-q_{2}}{{\alpha}}=\frac{49\alpha-36+\sqrt{144-264\alpha+289\alpha^2}}
{\alpha+12+\sqrt{144-264\alpha+289\alpha^2}}.$$ We can check that
$r$ is an increasing function with $\alpha\in(9/5,\infty]$, and
satisfies $$2<r\leq 11/3.$$ Therefore, for such $q_2,r, \alpha$, we
can apply Lemma 2.3 in \cite{[201]}, and  get
\begin{equation} \label{d14}\begin{array}{ll}
&\|\nabla
u\|_{L^2}^{2}+\displaystyle\frac{\nu}{2}\displaystyle\int_{0}^{t}\|\Delta
u\|_{L^2}^{2}d\tau\\
&\ \ \ \ \ \  \leq\|\nabla u(0)\|_{L^2}^{2}+\displaystyle
C\int_{0}^{t}\| u_{3}\|_{L^{q_2}}^{\frac{2(r-1)}{r-2}}\|
\partial_{3}u_{3}\|_{L^\alpha}^{\frac{2}{r-2}}\|\nabla
u\|_{L^2}^{2}d\tau\\
&\hspace{1.3cm} +C\displaystyle\int_{0}^{t}\|
u_{3}\|_{L^{q_2}}^{\frac{8(r-1)}{3(r-2)}}\|
\partial_{3}u_{3}\|_{L^\alpha}^{\frac{8}{3(r-2)}}\|\nabla
u\|_{L^2}^{2}d\tau+C.
\end{array}\end{equation}
Moreover, by the definition of $\mu$ and $r$, we have
\begin{equation} \label{d14}\begin{array}{ll}
\displaystyle\frac{8}{3(r-2)}-\frac{2\mu}{3-\mu}
&\displaystyle =\frac{8(\mu+\alpha)}{3(2\mu\alpha-3\mu-\alpha)}-\frac{2\mu}{3-\mu}\vspace{2mm}\\
&\displaystyle=\frac{2[12\mu+5\mu^2+(-6\mu^2-\mu+12)\alpha]}{3(2\mu\alpha-3\mu-\alpha)(3-\mu)}\\
&\displaystyle\equiv0.
\end{array}\end{equation}
 Combining \eqref{d13} and
\eqref{d14}, and the fact $\frac{2}{r-2}<\frac{8}{3(r-2)},$ we have
$$\begin{array}{ll}
\|\nabla
u\|_{L^2}^{2}+\displaystyle{\nu}\displaystyle\int_{0}^{t}\|\Delta
u\|_{L^2}^{2}d\tau\leq\displaystyle\displaystyle C V_{1}(t)+\|\nabla u(0)\|_{L^2}^{2}+C,\\
\end{array}
$$
and end the proof for $ \alpha\in\left(\frac{9}{5},\infty\right]$ by
using  Gronwall's inequality.\\
\noindent\textbf{Proof of Theorem \ref{t1.6}} The proof of this
theorem is heavily rely on the Lemma 2.3 in \cite{[201]}.\\
$\bullet \ i\neq j$\par For $(a)$, without loss of generality, we
assume $i=1,j=3$. The case of $\epsilon=1$ has been proved in
Theorem 1.1 in \cite{[201]}. For each  $1<\epsilon<3/2$ and
$$\alpha_1\in(\frac{3}{2-\epsilon},\frac{3}{3-2\epsilon}],$$ we take
$$
q=\frac{(12-4\epsilon)\alpha_1-12}{3(\alpha_1-1)},
$$
 $q$ is an increasing function with the variable
$\alpha$, and
$$
\frac{4}{1+\epsilon}<q\leq 2.
$$
By the initial data, using the similar argument to the proof of Theorem \ref{t1.4},
we see that
\begin{equation}\label{ttd15}
u_3\in L^{\infty}(0, T; L^{q}(\mathbb{R}^{3}))\ \mbox{with}\
1<q\leq2.
\end{equation}
Now, let $r=\frac{(q+1){\alpha_1}-q}{{\alpha_1}}$, then
$$
r=\frac{15\alpha_1-4\epsilon\alpha_1-12}{3\alpha_1}.
$$
We claim that $r>7/3$. In fact, it follows from
$$
\alpha_1>\frac{3}{2-\epsilon}\Longrightarrow
4\left(\frac{(2-\epsilon)\alpha_1}{3}-1\right)=
(q-\frac{4}{3})\alpha_1-q>0 \Longleftrightarrow r>\frac{7}{3}.$$
 Therefore, apply Lemma 2.3 in
\cite{[201]}, we get
\begin{equation}\label{t10} \begin{array}{ll}
&\|\nabla
u\|_{L^2}^{2}+\displaystyle{\nu}\displaystyle\int_{0}^{t}\|\Delta
u\|_{L^2}^{2}d\tau\vspace{2mm}\\&\hspace{1cm}\leq\displaystyle\displaystyle
C\int_{0}^{t}\|u_{3}\|_{L^q}^{\frac{8(r-1)}{3r-7}}\|
\partial_{1}u_{3}\|_{L^{\alpha_1}}^{\frac{8}{3r-7}}\|\nabla
u\|_{L^2}^{2}d\tau+\|\nabla u(0)\|_{L^2}^{2}\vspace{2mm}\\
&\hspace{1cm}\hspace{0.5cm}+\displaystyle
C\int_{0}^{t}\|u_{3}\|_{L^q}^{\frac{2(r-1)}{r-2}}\|
\partial_{1}u_{3}\|_{L^{\alpha_1}}^{\frac{2}{r-2}}\|\nabla
u\|_{L^2}^{2}d\tau+C.
\end{array}\end{equation}
From \eqref{ttd15}, it follows that
\begin{equation}\label{tt10} \begin{array}{ll}
&\|\nabla
u\|_{L^2}^{2}+\displaystyle{\nu}\displaystyle\int_{0}^{t}\|\Delta
u\|_{L^2}^{2}d\tau\vspace{2mm}\\&\hspace{1cm}\leq\displaystyle\displaystyle
C\int_{0}^{t}\|
\partial_{1}u_{3}\|_{L^{\alpha_1}}^{\frac{8}{3r-7}}\|\nabla
u\|_{L^2}^{2}d\tau+\|\nabla u(0)\|_{L^2}^{2}\vspace{2mm}\\
&\hspace{1cm}\hspace{0.5cm}+\displaystyle C\int_{0}^{t}\|
\partial_{1}u_{3}\|_{L^{\alpha_1}}^{\frac{2}{r-2}}\|\nabla
u\|_{L^2}^{2}d\tau+C,
\end{array}\end{equation}
it is obvious that $\frac{8}{3r-7}>\frac{2}{r-2}.$ Therefore, by
H$\ddot{\mbox{o}}$lder's inequality one has
\begin{equation}\label{ttt10} \begin{array}{ll}
&\|\nabla
u\|_{L^2}^{2}+\displaystyle{\nu}\displaystyle\int_{0}^{t}\|\Delta
u\|_{L^2}^{2}d\tau\leq\displaystyle\displaystyle C\int_{0}^{t}\|
\partial_{1}u_{3}\|_{L^{\alpha_1}}^{\frac{8}{3r-7}}\|\nabla
u\|_{L^2}^{2}d\tau+\|\nabla u(0)\|_{L^2}^{2}+C.
\end{array}\end{equation}
We see that
$$
\frac{3}{\alpha_1}+\frac{3r-7}{4}=\frac{3}{\alpha_1}+\frac{2\alpha_1-\epsilon\alpha_1-3}{\alpha_1}=2-\epsilon.$$
 By using  Gronwall's inequality and condition \eqref{tt15}, we prove
$(a)$.\\
$\bullet \ i=j$\par Without loss of generality, we assume $i=j=3$.
Firstly, we consider $s=1/2$.  
Let
\begin{equation}\label{sw1}
\frac{1}{\mu}+\frac{1}{2}=\frac{2}{q}\ \mbox{with}\
2<q\leq\frac{12}{5},
\end{equation}
from \eqref{sw1}, we see that $2<\mu\leq 3$. For above $q$, we prove
the following fact
\begin{equation}\label{sw2}
u_3\in L^{\infty}(0, T; L^{q}(\mathbb{R}^{3})).
\end{equation}
In fact, for such $\mu, q$, we can apply the same method to
\eqref{e24} (or see \eqref{d12}) and combine  Gagliardo-Nirenberg
and H$\ddot{\mbox{o}}$lder's inequalities to get
\begin{equation}\label{sw3}
\begin{array}{ll}
\displaystyle
\frac{1}{2}\frac{d}{dt}\|u_{3}\|_{L^q}^{2}\displaystyle
\leq\displaystyle C\| u\|_{L^2}^{\frac{3-\mu}{\mu}}\|\nabla
 u\|_{L^2}^{\frac{3(\mu-1)}{\mu}}\|\partial_{3}u_{3}\|_{L^2}.
\end{array}
\end{equation}
Integrating \eqref{sw3} on time, applying energy inequality and the
condition \eqref{tts15}, one has
\begin{equation}\label{sw4}
\begin{array}{ll}
\displaystyle \|u_{3}\|_{L^q}^{2}\displaystyle &\leq\displaystyle
\|u_3(0)\|_{L^q}+C\int_0^t\| u\|_{L^2}^{\frac{3-\mu}{\mu}}\|\nabla
 u\|_{L^2}^{\frac{3(\mu-1)}{\mu}}\|\partial_{3}u_{3}\|_{L^2}d\tau\vspace{1mm}\\
 \displaystyle &\leq\displaystyle\|u_3(0)\|_{L^q}+C\int_0^t\| u\|_{L^2}^{2}d\tau+C\int_0^t\|\nabla
 u\|_{L^2}^{2}d\tau\vspace{1mm}\\
 \displaystyle &\leq\displaystyle\|u_3(0)\|_{L^q}+C(T).
\end{array}
\end{equation}
This proves \eqref{sw2}. Let
$$
r=\frac{q+2}{2},
$$
then $2<r\leq11/5,$ and $\frac{q}{q-r+1}=2$.  Taking the inner
product of the equation \eqref{a} with $-\Delta_{h}u$ in $L^{2}$,
applying H$\ddot{\mbox{o}}$lder's inequality several times, we
obtain
\begin{equation}\label{sw5}
\begin{array}{ll}
 &\displaystyle \frac{1}{2}\frac{d}{dt}\|\nabla_{h}u\|_{L^2}^{2}+\nu\|\nabla_{h}\nabla u\|_{L^2}^{2}\\
 &\ \ \ \ =\displaystyle
 \int_{\mathbb{R}^{3}}[(u\cdot \nabla)u]\Delta_{h}u dx\displaystyle \\
 \displaystyle &\ \ \ \ \leq\displaystyle
 C\int_{\mathbb{R}^{3}}|u_{3}||\nabla u||\nabla_{h}\nabla u|dx \ \ (\mbox{see}\  \mbox{\cite{[2]}})\displaystyle \\
 \displaystyle &\ \ \ \ \leq\displaystyle
 C\int_{\mathbb{R}^{2}}\max_{x_{3}}|u_{3}|(\int_{\mathbb{R}}|\nabla u|^{2}dx_{3})^{\frac{1}{2}}(\int_{\mathbb{R}}|\nabla_{h}\nabla u|^{2}dx_{3})^{\frac{1}{2}} dx_{1}dx_{2}\
 \displaystyle \vspace{2mm}\\
 \displaystyle &\ \ \ \ \leq\displaystyle
 C[\int_{\mathbb{R}^{2}}(\max_{x_{3}}|u_{3}|)^{r}dx_{1}dx_{2}]^{\frac{1}{r}}[\int_{\mathbb{R}^{2}}(\int_{\mathbb{R}}|\nabla
 u|^{2}dx_{3})^{\frac{r}{r-2}}dx_{1}dx_{2}]^{\frac{r-2}{2r}}\\
 \displaystyle& \ \ \ \ \ \ \times \displaystyle[\int_{\mathbb{R}^{3}}|\nabla_{h}\nabla
 u|^{2}dx_{1}dx_{2}dx_{3}]^{\frac{1}{2}}\\
 \displaystyle \vspace{2mm}
 \displaystyle &\ \ \ \ \leq\displaystyle
 C[\int_{\mathbb{R}^{3}}|u_{3}|^{r-1}|\partial_{3}u_{3}|dx_{1}dx_{2}dx_{3}]^{\frac{1}{r}}\|\nabla_{h}\nabla u\|_{2}
  \displaystyle \\
 \displaystyle& \ \ \ \ \ \ \times \displaystyle[\int_{\mathbb{R}}(\int_{\mathbb{R}^{2}}|\nabla
 u|^{\frac{2r}{r-2}}dx_{1}dx_{2})^{\frac{r-2}{r}}dx_{3}]^{\frac{1}{2}}
 \displaystyle \vspace{2mm}\\
\displaystyle &\ \ \ \ \ \leq \displaystyle C\displaystyle\|
u_{3}\|_{L^q}^{\frac{r-1}{r}}\|
\partial_{3}u_{3}\|_{L^{\frac{q}{q-r+1}}}^{\frac{1}{r}}\|\nabla u\|_{L^2}^{\frac{r-2}{r}}
 \|\partial_{2}\nabla u\|_{L^2}^{\frac{1}{r}} \|\partial_{3}\nabla u\|_{L^2}^{\frac{1}{r}}
\|\nabla_{h}\nabla u\|_{L^2}.\vspace{2mm}\\
\displaystyle &\ \ \ \ \ =\displaystyle C\displaystyle\|
u_{3}\|_{L^q}^{\frac{r-1}{r}}\|
\partial_{3}u_{3}\|_{L^2}^{\frac{1}{r}}\|\nabla u\|_{L^2}^{\frac{r-2}{r}}
 \|\partial_{2}\nabla u\|_{L^2}^{\frac{1}{r}} \|\partial_{3}\nabla u\|_{L^2}^{\frac{1}{r}}
\|\nabla_{h}\nabla u\|_{L^2}.
\end{array}\end{equation}
Applying Young's inequality to \eqref{sw5}, we have
\begin{equation}\label{sw6}
\begin{array}{ll}
 &\displaystyle \frac{1}{2}\frac{d}{dt}\|\nabla_{h}u\|_{L^2}^{2}+\nu\|\nabla_{h}\nabla u\|_{L^2}^{2}\\
\displaystyle &\ \ \ \ \ \leq \displaystyle C\displaystyle\|
u_{3}\|_{L^q}^{\frac{2(r-1)}{r-2}}\|
\partial_{3}u_{3}\|_{L^2}^{\frac{2}{r-2}}\|\nabla u\|_{L^2}^{2}+\frac{\nu}{2}\|\nabla_{h}\nabla u\|_{L^2}^{2}.
\end{array}\end{equation}
After integrating \eqref{sw6} on time, combining the energy
inequality and \eqref{sw2}, as well as the condition \eqref{tts15},
one has
\begin{equation}\label{sw7}
\begin{array}{ll}
 &\displaystyle \|\nabla_{h}u\|_{L^2}^{2}+\nu\int_0^t\|\nabla_{h}\nabla u\|_{L^2}^{2}d\tau\\
\displaystyle &\ \ \ \ \ \ \ \ \ \ \leq \displaystyle
C\int_0^t\displaystyle\| u_{3}\|_{L^q}^{\frac{2(r-1)}{r-2}}\|
\partial_{3}u_{3}\|_{L^2}^{\frac{2}{r-2}}\|\nabla u\|_{L^2}^{2}d\tau+\|\nabla_{h}
u(0)\|_{L^2}^{2}.\vspace{1mm}\\
\displaystyle &\ \ \ \ \ \ \ \ \ \  \leq \displaystyle
C(T)\int_0^t\displaystyle\|\nabla u\|_{L^2}^{2}d\tau+\|\nabla_{h}
u(0)\|_{L^2}^{2}.\vspace{2mm}\\
&\ \ \ \ \ \ \ \ \ \  \leq \displaystyle C(T)+\|\nabla_{h}
u(0)\|_{L^2}^{2}.
\end{array}\end{equation}
On the other hand, we have (see  the proof of Theorem 1.3 in
\cite{[t201]} for detail)
\begin{equation} \label{sw9}\begin{array}{ll}\displaystyle \frac{1}{2}\frac{d}{dt}\|\nabla
u\|_{L^2}^{2}+\nu\|\Delta u\|_{L^2}^{2}\leq \displaystyle
C\|\nabla_{h}u\|_{L^2}\|\nabla
u\|_{L^2}^{\frac{1}{2}}\|\nabla_{h}\nabla u\|_{L^2}\|\Delta
u\|_{L^2}^{\frac{1}{2}}.
\end{array}\end{equation}
After integrating, and  using \eqref{sw7} and energy inequality, we
obtain
\begin{equation} \label{sw8}\begin{array}{ll}
&\|\nabla
u\|_{L^2}^{2}+\displaystyle2\nu\displaystyle\int_{0}^{t}\|\Delta
u\|_{L^2}^{2}d\tau\\
\displaystyle &\hspace{0.3cm} \leq\|\nabla
u(0)\|_{L^2}^{2}+\left(\sup_{0\leq s\leq t}\|\nabla_{h}
u\|_{L^2}\right)\left(\displaystyle\int_{0}^{t}\|\nabla
u\|_{L^2}^{2}d\tau \right)^{\frac{1}{4}}\\
&\ \ \ \ \ \times\left(\displaystyle\int_{0}^{t}\|\nabla_{h}\nabla
u\|_{L^2}^{2}d\tau\right)^{\frac{1}{2}}\left(\displaystyle\int_{0}^{t}\|\Delta
u\|_{L^2}^{2}d\tau\right)^{\frac{1}{4}}\\
\displaystyle &\hspace{0.3cm} \leq\|\nabla
u(0)\|_{L^2}^{2}+C(T)\left(\displaystyle\int_{0}^{t}\|\Delta
u\|_{L^2}^{2}d\tau\right)^{\frac{1}{4}}.
\end{array}\end{equation}
By Young inequality, we get the $ H^1$ norm of the strong solution
$u$ is bounded on the maximal interval of existence $(0, T^{*})$.
This completes the proof of the case of $\epsilon=1/2$.
\par For each
$$\alpha_2\in(\frac{3}{2-\epsilon},\frac{6}{5-4\epsilon}],\ \mbox{with} \ 1/2<\epsilon<5/4,$$ we take
$$
q=\frac{(11-4\epsilon)\alpha_2-12}{3(\alpha_2-1)}.
$$
It is easy to see that $q$ is an increasing function of $\alpha_2$
and
$$
\frac{3}{1+\epsilon}<q\leq2\ (\frac{1}{2}<\epsilon<\frac{5}{4}).
$$
By using the initial data, as before (see the proof in Theorem
\ref{t1.4} for detail), we have
\begin{equation}\label{tptd15}
u_3\in L^{\infty}(0, T; L^{q}(\mathbb{R}^{3}))\ \mbox{with}\
1<q\leq2.
\end{equation}
Let $r=\frac{(q+1){\alpha_2}-q}{{\alpha_2}}$, then one has
$$
r=\frac{14\alpha_2-4\epsilon\alpha_2-12}{3\alpha_2},
$$
and $r>2$. In fact, it follows from
$$
r-2=\frac{4(2\alpha_2-\epsilon\alpha_2-3)}{3\alpha_2}>0\
(\alpha_2>\frac{3}{2-\epsilon}).
$$
Therefore, we can apply Lemma 2.3 in \cite{[201]}, and get
\begin{equation} \label{tc95}\begin{array}{ll}
&\|\nabla
u\|_{L^2}^{2}+\displaystyle\frac{\nu}{2}\displaystyle\int_{0}^{t}\|\Delta
u\|_{L^2}^{2}d\tau\\
&\ \ \ \ \ \  \leq\|\nabla u(0)\|_{L^2}^{2}+\displaystyle
C\int_{0}^{t}\| u_{3}\|_{L^q}^{\frac{2(r-1)}{r-2}}\|
\partial_{3}u_{3}\|_{L^{\alpha_2}}^{\frac{2}{r-2}}\|\nabla
u\|_{L^2}^{2}d\tau\\
&\hspace{1.3cm} +C\displaystyle\int_{0}^{t}\|
u_{3}\|_{L^q}^{\frac{8(r-1)}{3(r-2)}}\|
\partial_{3}u_{3}\|_{L^{\alpha_2}}^{\frac{8}{3(r-2)}}\|\nabla
u\|_{L^2}^{2}d\tau+C.
\end{array}\end{equation}
From \eqref{tptd15}, it follows that
\begin{equation}\label{tpt10} \begin{array}{ll}
&\|\nabla
u\|_{L^2}^{2}+\displaystyle\frac{\nu}{2}\displaystyle\int_{0}^{t}\|\Delta
u\|_{L^2}^{2}d\tau\\
&\ \ \ \ \ \  \leq\|\nabla u(0)\|_{L^2}^{2}+\displaystyle
C\int_{0}^{t}\|
\partial_{3}u_{3}\|_{L^{\alpha_2}}^{\frac{2}{r-2}}\|\nabla
u\|_{L^2}^{2}d\tau\\
&\hspace{1.3cm} +C\displaystyle\int_{0}^{t}\|
\partial_{3}u_{3}\|_{L^{\alpha_2}}^{\frac{8}{3(r-2)}}\|\nabla
u\|_{L^2}^{2}d\tau+C,
\end{array}\end{equation}
it is obvious that $r>2$ implies $\frac{8}{3r-6}>\frac{2}{r-2}.$
Therefore, by H$\ddot{\mbox{o}}$lder's inequality one has
\begin{equation}\label{ttt10} \begin{array}{ll}
&\|\nabla
u\|_{L^2}^{2}+\displaystyle{\nu}\displaystyle\int_{0}^{t}\|\Delta
u\|_{L^2}^{2}d\tau\leq\displaystyle\displaystyle C\int_{0}^{t}\|
\partial_{3}u_{3}\|_{L^{\alpha_2}}^{\frac{8}{3(r-2)}}\|\nabla
u\|_{L^2}^{2}d\tau+\|\nabla u(0)\|_{L^2}^{2}+C.
\end{array}\end{equation}
We see that
$$
\frac{3}{\alpha_2}+\frac{3r-6}{4}=\frac{3}{\alpha_2}+\frac{2\alpha_2-\epsilon\alpha_2-3}{\alpha_2}=2-\epsilon.$$
 By
using  Gronwall's inequality and condition \eqref{tts15}, we prove
$(b)$.\\
\begin{rem}
In the proof of $(b)$ of this theorem, the result of case
$\epsilon=1/2$ is actually obtained in Theorem \ref{t1.5}, in which
we note that $\beta_2$ is not necessary to be infinity when
$\alpha_2=2$. However, when $\beta_2=\infty$ we have a clear proof,
and we have shown above.
\end{rem}
\noindent\textbf{Proof of Theorem \ref{t1.7}} The method of the
proof of $(a)$ is quite similar to the Theorem 1.2 (i) in
\cite{[201]}, therefore, we only the give the outline of the proof.
We also assume $i=1,\ j=3$. For every
$$
\alpha_1\in[\frac{3}{3-2\epsilon},\frac{3(11-2\epsilon)}{2\epsilon^2-26\epsilon+33}],
\ \mbox{with}\ \  1<\epsilon\leq21/16,
$$
 we set
\begin{equation}\label{jdsw2}\left\{\begin{array}{l}
\displaystyle q=\frac{(12-4\epsilon)\alpha_1-12}{3(\alpha_1-1)},\vspace{2mm}\\
\displaystyle \sigma=\frac{3\alpha_1q}{6\alpha_1-q},\vspace{2mm}\\
\displaystyle r=\frac{(q+1)\alpha-q}{\alpha},
\end{array}
\right.
\end{equation}
then we have $2<q<6$, and
$$
\sigma=\frac{(18-6\epsilon)\alpha_1^2-18\alpha_1}{9\alpha_1^2-(15-2\epsilon)\alpha_1+6},
r=\frac{(15-4\epsilon)\alpha_1-12}{3\alpha_1},
$$
$$
\frac{1}{\sigma}+\frac{q-2}{q}+\frac{1}{3{\alpha_1}}=1,\ 1<\sigma<
\frac{9}{8},
$$
moreover, we also have  $r>\frac{7}{3},$ and
\begin{equation}\label{jws2}
\frac{8}{3r-7}>\frac{2\sigma}{9-8\sigma}\ \mbox{with}\
\alpha_1\in[\frac{3}{3-2\epsilon},\frac{3(11-2\epsilon)}{2\epsilon^2-26\epsilon+33}],
\ \mbox{and}\ \  1<\epsilon\leq21/16.
\end{equation}
For the rest, we will use these parameters $q,\sigma, r$ to  get the
desired result. As give the same process to \cite{[201]}, we will
use Lemma 2.1 of \cite{[201]} to estimate $u_3$ with the parameters
$q,\sigma$, and then by the $r>7/3$ to get Lemma 2.3 (i) in
\cite{[201]}. Finally, by \eqref{jws2} and condition \eqref{jt15} to
finish the proof.
\par Now, we pay attention to $(b)$, and assume $i=j=3.$ The proof of this part is to imitate
  the proof of  Theorem \ref{t1.5}, and for every
$$\alpha_2\in[\frac{6}{5-4\epsilon},\frac{18-2\epsilon}{(4\epsilon-3)(\epsilon-5)}]\ \mbox{with}\ \frac{1}{2}<\epsilon
\leq\frac{3}{4},$$ we set
\begin{equation}\label{dsw2}\left\{\begin{array}{l}
\displaystyle q=\frac{(11-4\epsilon)\alpha_2-12}{3(\alpha_2-1)},\vspace{2mm}\\
\displaystyle \mu=\frac{\alpha_2 q}{2\alpha_2-q},\vspace{2mm}\\
\displaystyle r=\frac{(q+1)\alpha_2-q}{\alpha_2}.
\end{array}
\right.
\end{equation}
From the definition, we can check that $2<q<6$, and
$$
\mu=\frac{(11-4\epsilon)\alpha_2^2-12\alpha_2}{6\alpha_2^2-(17-4\epsilon)\alpha_2+12}
,r=\frac{2(7-2\epsilon)\alpha_2-12}{3\alpha_2},
$$
$$
\frac{1}{\mu}+\frac{1}{\alpha_2}+\frac{q-2}{q}=1, \ 1<\mu<3,\ 2<q<6,
$$
and  $r>2$, moreover, we also note the index satisfies
\begin{equation}\label{jjjdsw2}
\frac{8}{3(r-2)}>\frac{2\mu}{3-\mu}\ \ \mbox{with}\
\alpha_2\in[\frac{6}{5-4\epsilon},\frac{18-2\epsilon}{(4\epsilon-3)(\epsilon-5)}]\
\mbox{and}\ \frac{1}{2}<\epsilon\leq\frac{3}{4}. \end{equation}
Finally, we will use these parameters $q,\mu, r$ to give the same
process as the proof of Theorem \ref{t1.5} to get this result. We
will bounds the $\|u_3\|_{L^q}$ and $\|\nabla u\|_{L^{2}}$ in turn
by  \eqref{jjjdsw2} and condition \eqref{tsw15}. We complete the
proof.

\section*{Acknowledgement}
The second author would like to thank Dr.Ting Zhang for his helpful
suggestions. This work is supported partially by
NSFC 11271322, 10931007, and Zhejiang NSF of China Z6100217, and LR12A01002.

}
\end{document}